


%
%
%
%

\baselineskip=.35in

\let\bold\bf

\def\({\big(}
\def\){\big)}
\def\[{\big[}
\def\]{\big]}

\def\0{\emptyset}

\def\nfork{\setbox0\hbox{$\bigcup$}%
\setbox1=\hbox to \wd0{\hfil$|$\hfil}%
\wd1=0cm\relax\box1\box0}

\def\beth{%
\hbox{\vrule width1pt height0ptdepth0pt
\vrule width 4pt depth-5pt height6pt \vrule width1pt
height6pt depth 0pt \kern-6pt\vrule height1pt width 7pt}}


\newfam\frakturfam
\newfam\msyfam
\newfam\msxfam

\def\amsfonts{%
\font\teneuf=eufm10
\font\seveneuf=eufm7
\font\fiveeuf=eufm5
\textfont\frakturfam\teneuf
\scriptfont\frakturfam\seveneuf
\scriptscriptfont\frakturfam\fiveeuf
\font\msbm=msbm10
\font\sevenmsbm=msbm7
\font\fivemsbm=msbm5
\textfont\msyfam=\msbm
\scriptfont\msyfam=\sevenmsbm
\scriptscriptfont\msyfam=\fivemsbm
\font\msam=msam10
\font\sevenmsam=msam7
\font\fivemsam=msam5
\textfont\msxfam=\msam
\scriptfont\msxfam=\sevenmsam
\scriptscriptfont\msxfam=\fivemsam
}

\def\frak{\fam\frakturfam\teneuf}
\let\goth\frak

\def\Bbb{\relax\ifmmode\fam\msyfam\else\message{Bbb not allowed in
text}\fi}
\def\Bbb#1{\relax\ifmmode\fam\msyfam#1\else$\Bbb #1$\fi}

\def\hexnumberat#1{\ifcase#1 0\or 1\or 2\or 3\or 4\or 5\or 6\or 7\or 8\or
 9\or A\or B\or C\or D\or E\or F\fi}
\edef\msxm{\hexnumberat\msxfam}
\edef\msym{\hexnumberat\msyfam}

\mathchardef\upharpoonright="3\msxm16
\let\restriction\upharpoonright

\mathchardef\notvDash="3\msxm32
\mathchardef\Vdash="3\msxm0D
\mathchardef\beth="0\msym69
\mathchardef\notVdash="3\msym31
\mathchardef\notvdash="3\msym30
\mathchardef\triangleleft="3\msxm43
\mathchardef\trianglelefteq="3\msxm45
\mathchardef\triangleright="3\msxm42
\mathchardef\trianglerighteq="3\msxm44

\def\lessdot{\mathrel{\mathord{<}\!\!%
  \raise 0.8 pt\hbox{$\scriptstyle\circ$}}}
\def\lesseqdot{\mathrel{\mathord{\le}\!\!%
  \raise 0.8 pt\hbox{$\scriptstyle\circ$}}}

\def\monthname#1{\ifcase#1\or
Jan\or
Feb\or
Mar\or
Apr\or
May\or
Jun\or
Jul\or
Aug\or
Sep\or
Oct\or
Nov\or
Dec\fi}
\def\date{\monthname\month\ \number\day, \number\year}

\footline={\tenrm \date\ (\jobname)\hfill \folio\hfill \ \ \ \ }

\parindent0cm
\advance\hsize1cm

\def\draft{\magnification\magstep2
\baselineskip\vsize
\advance\baselineskip-\topskip
\divide\baselineskip18
\amsfonts
\mylines}

\def\final{\magstep1
\amsfonts
\mylines}



\def\XX{\mathrel{\hbox{\msam\char"43}}}
\def\notXX{\;\>\rlap/\kern-5pt\XX}

\def\lk{\langle}
\def\rk{\rangle}

\catcode`\"\active
\def"{\ifleftquote\lq\global\leftquotefalse\else\rq\global\leftquotetrue\fi}
\newif\ifleftquote
\leftquotetrue

\def\lq{\hbox{\rm``}}
\def\rq{\hbox{\rm''}}

\def\qed{\hfill\hbox{\msam\char'003}}
\def\and{{\, \&\,}}

\def~{\hbox{$\sim$}}

\def\page{\vfill\eject}


{\catcode`\ \active
\global\def {\leavevmode\space}}
\everymath={\catcode`\ 10 }

\def\endgrafleavevmode{\leavevmode\endgraf}

\def\mylines{
\let\par\endgrafleavevmode%
\obeylines%
\catcode`\ \active}
\def\texlines{%
\let\par\endgraf%
\relax%
\catcode`\ 10\relax%
\catcode`\^^M5\relax}

\output={\let\par\endgraf\plainoutput}
\catcode`\#12






%
%



\def\bethA{\hbox{\vrule width1pt
height0ptdepth0pt
\vrule width 4pt depth-5pt height6pt \vrule width1pt
height6pt depth 0pt \kern-6pt\vrule height1pt width 7pt}}
\def\beth{\Bbb{i}}

\def\st{ such that }

\def\a{\alpha}
\def\A{\Rightarrow}



\def\h{\eta}

\def\L{\Lambda}

\def\n{\nu}

\def\p{ \pi }

\def\R{\varrho}
\def\r{\rho}

\def\S{\Sigma}
\def\t{\tau}

\def\w{\omega}

\def\phi{\varphi}

\def\rt{\upharpoonright}

\def\p0{\emptyset}
\def\pa{\forall}
\def\pb{\setminus}
\def\pc{\wedge}

\def\pE{\prec}
\def\pf{\Vdash}

\def\pl{\vert}

\def\pr{\upharpoonright}

\def\ps{\subseteq}

\def\pu{\cup}

\def\pv{\cap}

\def\py{\exists}

\ifx\mathcal\undefinedcs
 \def\mathcal{\cal}\fi
\ifx\cal\undefinedcs\def\cal{\Cal}\fi

\def\sH{\mathcal{H}}

\ifx\mathfrak\undefinedcs\def\mathfrak{\frak}\fi

\def\gc{{\mathfrak{c}}}

\ifx\mathbf\undefinedcs\def\mathbf{\bf}\fi

\def\bp{\mathbf{p}}

\def\bx{\mathbf{x}}

\def\xB{{\Bbb{B}}}

\def\xE{{\Bbb{E}}}

\def\xN{\Bbb{N}}

\def\xP{{\Bbb{P}}}

\def\xQ{{\Bbb{Q}}}

\def\xR{{\Bbb{R}}}

\def\xZ{\Bbb{Z}}


\def\0{emptyset}

\def\lk{\langle}
\def\rk{\rangle}

\def\({\big(}
\def\){\big)}
\def\[{\big[}
\def\]{\big]}

\def\0{\emptyset}

\def\lessdot{\mathrel{\mathord{<}\!\!%
  \raise 0.8 pt\hbox{$\scriptstyle\circ$}}}
\def\lesseqdot{\mathrel{\mathord{\le}\!\!%
  \raise 0.8 pt\hbox{$\scriptstyle\circ$}}}

\def\nfork{\setbox0\hbox{$\bigcup$}%
\setbox1=\hbox to \wd0{\hfil$|$\hfil}%
\wd1=0cm\relax\box1\box0}


 \final
 \texlines




{\bold A Partition Theorem }

 Saharon Shelah

\endgraf
Institute of Math, The Hebrew University, Jerusalem
Israel
\endgraf
Department of Math., Rutgers University, New
Brunswick NJ USA

footnote Presented in
the third Turan lecture of
the author
 in Hungary, Feb 1998;
 Publ.  Number 679. Done 1-2/98. Partially supported
by the  Partially supported by
the  United States Israel Binational Science Foundation.

{\bold Abstract}

We deal with some relatives of the
Hales Jewett theorem with primitive recursive
bounds.

\

{\bold Anotated Content}

0 Introduction

1  Basic Definitions

2 Proof of the Partition Theorem with a bound

3 : Higher Dimension Theorems

4 The Main Theorem

\


\

Classification: Primary 05A99, Secondary 15A03

 {\bold  0: Introduction }

 We prove the following: there is a primitive
recursive function $f_{-}^*(-,- )$ , in the three
variables, such that:
for every natural numbers
$t, n>0$, and $c$,
  for any natural number $k \ge f^*_t (n,c)$
 the following holds. Assume  $\Lambda$
is an alphabet
 with $n>0$ letters,
 $M$ is the family of non empty
 subsets of $\{1,\dots, k\}$ with $\le t$ members
 and $V$ is the  set of functions from $M$
 to $\Lambda$ and lastly
 $d$ is a $c-$colouring of $V$ (i.e.
 a function with domain
 $V$ and range with at most $c$ members). Then
  there is a $d-$monochromatic $V-$line, which means that
 there are  $w \subseteq \{ 1,\dots,k\}$,
with at least
 $t$ members and function $\rho$ from
 $\{ u \in M: u$ not a subset of $w\}$  to $\Lambda$
 such that letting $L = \{ \eta \in V : \eta$ extend $\rho$
and  for each $s=1,\dots,t$
 it is constant
 on $\{ u \in M: u \subseteq w$ has $s$ members $\}\}$,
we have $d \restriction L$ is constant
(for $t=1$ those are the Hales Jewett numbers).

A second theorem relates to the first
just as the affine Ramsey theorem of
Graham, Leob and Rothschild (which continue
the n-parameter Ramsey theorem of
Graham and  Rothschild), relates
to the Hales Jewett theorem. We also
note an infinitary related theorem
parallel to the Galvin Prikry theorem
and the Carlson Sympson theorem.

Let us review history and background,
not repeating [GRS 80]. In the late seventies,
Furstenberg and Sarakozy independently
prove that if $\bp (x)$ is a polynomial
in $\xZ [\bx] $ satisfying $p(0)=0$
and $A \ps  \xN$ is a set of positive density
then for some $a,b \in A$ and  $n\in \xN$
we have $a-b = {\bp } (n)$. Bergelson and Leibman
 [BL96] continuing Furstenberg [Fu] prove
(this is a special of a density theorem like
Szemeredi ):
if $r,k,t,m$ are natural numbers,
${\bp }_{\ell,s} (x)$ for $\ell = 1, \dots , k$
and $s=1,\dots ,t$ are polynomials with rational
coefficients, taking integer values
at integers, and vectors
${\bar v}_1,\dots , {\bar v}_t \in {}^m\xZ$
and any $r-$colouring of ${}^m\xZ$ there
are $ \bar a \in {}m\xZ$ and $n \in {\xZ} (n \neq 0)$
\st the set
 $S(\bar a , n ) = \{ {\bar a} +
\S_{j=1,t} {\bp}_{i,j}  (n)  {\bar  v}_j :
i = 1, \dots , k \}$
is monochromatic.

Bergelson and Leibman [ BL 9x] prove a theorem,
"set polynomial extension", which is,
in different formulation, like the first
theorem describe above but without a bound
(i.e. the primitive recursiveness).
Their method is infinitary so does not seem to give
even the weak bound in 2.5 (one with triple
induction), and certainly does
not give primitive recursive bounds.

Naturally our proofs continue [Sh 329]. We thank
the referee for telling us on [BL 9x] and other
helpful comments.See a discussion of related problems in
[Sh 702].

 \

 0.1 { {\bold Notation}}:

 (a) We use $\Lambda$ for a finite alphabet, always non empty,
  members of which are denoted by $\alpha,\beta,\gamma.$

 (b) We use $M,N$   to denote structures which serve as index sets,
 so we call them index models. We use $\tau$ to denote vocabularies,
 , (see Definition 1.1), $F$ to denote function symbols.

 (c) We use $n,m,k,\ell, i,j,c,r,s,t$  to denote natural numbers,
 but usually $n$ is the number of  letters, i.e.
the number of members in an alphabet;
 $k$ the dimension of the index models
 and $c\ge 1$ the number of colours.

 (d) $|X|$ and also card$(X)$  denote
  the number of elements of the set $X.$

 (e)  We use $\eta,\nu,\rho$ to denote
members of spaces,
 we use $V,U$ to denote spaces
and $a,b$ to denote elements of $M,N$
 and $d$ to denote colourings, $p$ to denote the
 `type` of a point in a line and {\goth p}
 to denote
 type of a line or a space
 (see Definition 1.7(3)).
  We use $L$ to denote  (combinatorial)
lines, $S$ to denote (combinatorial)
 subspaces.

 (f) $A$ bar on a symbol, say $\bar x$ denote a
finite
 sequence of such objects,  of length
lg$(\bar x)$ the ???$i-$th object
 being $x_{i}$  (and of $\bar x_{m}$   or
$\bar x^{m}$  it is $x^{m}_{i} ).$

 \


 0.2 DEFINITION: (1) For $m \ge 1$,
 let {\Bbb E}$_{m}$  be  the minimal
 class of functions
 from natural numbers to natural numbers
 (with any number of places)
 closed under composition, which for $m=1$
 contains $0,1,x+1$
 and the projection functions,
and for $m>1$
contains
 any function which we get by inductive definition on
 functions from {\Bbb E}$_{m-1}$
 (see [Ro84], so {\Bbb E}$_3$ is the family of polynomials, {\Bbb E}$_4$
 contains the tower function and {\Bbb E}$_5$ contains
 the waw function  and $\cup_{m\ge 1} {\Bbb E}_{m}$ is the
 family of primitive recursive functions, and
 the `simplest` function not there is the Akerman function.)
 We allow an object like {$\bar \Lambda$} to be one of the arguments
 meaning a natural number coding of it (in the cases used this does not
matter).
 Abusing notation, we may say
"$f$ is in {\Bbb E}$_{n}$"
 instead of
"$f$ is bounded by a
function from {\Bbb E}$_{n}$", also
 writing $f_{\bar \Lambda}(-,\dots)$
we count $\bar \Lambda$
 as
one of
  the arguments.

 (2) We can define the Akerman
 function $A_{n}(m)$
by double induction
 (in as sense it is the simplest, smallest function which is not
 primitive recursive).

 \

 0.3 DEFINITION: (1) Let RAM$(t,\ell,c)$
be the
Ramsey number, i.e. the
  first $k$ such that  $k \rightarrow (t)^{\ell}_{c}$
 which mean  that if $A$ is a set with $k$ elements, and $d$ is
 a $c-$colouring of
$[A]^{\ell} =^{df} \{ B: B$ is a subset of $A$ with
$\ell$ elements$\}$, that is a function with this
domain and range
 of cardinality  $\le c$,
{\bold  then} for some
$A_1 \in [A]^{t}$ we have
$ d \restriction [A_1]^\ell$
 is constant.

 (2) Let HJ$(n,m,c)$ be the Hales Jewett number
 for getting a monochromatic
subspace of dimension $m$, when the colouring has $c$
 colours and
 for an alphabet with $n$ members (this is, by our
 subsequent definitions, $f^1(\bar \Lambda, m,c)$
 when $\tau(\bar \Lambda)= \{$id$\}$,
 and $\Lambda_{{\rm id}}$ has $n$ members, see Definition 1.9).

 \page

{\bold  Section 1 : Basic definitions }

We can look at Hales Jewett theorem
in geometric terms: $\xR$ is replaced by
$\L$; a finite alphabet, the $k-$dimensional
euclidean space ${\xR}^k$
 is replaced by ${}^{[1,k]} \L$
(or ${}^{[0,k)} \L$), essentially the set of
sequences of length $k$ of members of the alphabet
$\L$; a subspace is replaced by the set of solutions
$(x_1,\dots , x_k) \in {}^{[1,k]} \L$
of a family of linear equations, which here
means just $x_i = \a$ (where $\a \in \L , 1 \le i
\le k)$ or $x_i = x_j$. Here the basic set
$[1,k]$ is replaced by a structure $M$, a $\t-$fim.
Such  basic definitions are given in this section.

 We define a `space over an index model of  dimension $k$,
 over  an alphabet $\Lambda$ of size $n`$,
  lines and more.
 We then define the function $f^1$,  such that  for every  $n$, if $k$ is
 $f^1_{\tau}(n,c)$
 then for
 every colouring  of  the space by $\le c$ colours,
 there is a monochromatic line (in the appropriate interpretation.)
 Of course the use of id as a special function symbol
 is not really needed, also  we can waive the linear order on $P^{M}$,
 and the set of  automorphisms of the resulting
structure are natural for our
purpose,
 but not for the structures from 1.10(3); but at present those
 decisions does not matter.

 \

 1.1 DEFINITION:   (1) We call  $M$ a full index model  [fim
 or $\tau-$fim or fim for $\tau ]$ if:

 (a) the vocabulary $\tau = \tau_{M} = \tau(M)$ of $M$ includes a  unary
predicates $P$,
 a binary predicate $<$,
 and  finitely many function symbols   $F$,
 $F$ being arity$(F)-$place
 and no other symbols  (so $F$ vary over such function symbols).
 We may write arity$^{\tau} (F)$ for arity$(F).$
 We  usually treat $\tau$
 as the set of function symbols in $\tau.$

 (b)  the universe of $M$ is finite  (non empty of course).

 $(c) <^{M}$ is a linear order of $P^{M}$, so $x <^{M}  y$ implies $x,y
\in P^{M}.$

 $(d) F^{M}$ is a partial function  such that if $F^{M}(a_1,\dots,a_{r})$
 is well defined  (so $r =$ arity$(F) )$ then
  $a_1,\dots,a_{r} \in P^{M}$
 and the function is symmetric, i.e. does not depend
 on the order of the arguments, so if not said otherwise we assume
 $a_1 \le^{M}  a_2  \le^{M} \dots \le^{M} a_{r} .$

 (e) if $F^{M}_1(a_1,\dots,a_{r}) = F^{M}_2 (b_1,\dots,b_{t})$
 then $F_1 = F_2$ (hence $r=t )$  and $a_{\ell} = b_{\ell}$ for $\ell=
1,\dots,t$
 (under the convention from clause $(d))$
 and   every $b\in M \setminus  P^{M}$   has this form. So we  let
  base$_{M} (b) =^{df}   \{ a_1,\dots,a_{r}\}$
  and let base$_{\ell}(b)  =$ base$_{M,\ell} =^{df} a_{\ell}$
 where $b = F^{M} (a_1,\dots,a_{r})$
 (and $a_1 \le^{M}  a_2  \le^{M} \dots \le^{M} a_{r}$ of course)
 and   $F_{M,b}  =^{df}  F ;$ those are well defined by the
 demand above.

 $(f)$ $P^{M}$ is non empty and we call its cardinality dim$(M)$,
 the dimension of $M.$

 (g) id$^{M}$ is the identity function on $P^{M}$, so id is a
  unary function symbol of  $\tau$.

 (h) each $F^{M} (  a_1, \dots, a_{{\rm arity}(F)} )$  is well defined
 iff    $a_1, \dots, a_{{\rm arity}(F)}$  are from  $P^{M}$
 ( and the value does not depend on the order, of course)

 (2)  For $\tau$ as in part (1), let arity$(\tau)$  be Max$\{$arity$(F) :
F\in \tau\}$,
 so it is at least 1
 and let $\bar m[{\tau}] =^{df} \langle m^{\tau}_{t} :
t=1,\dots$,arity$(\tau) \rangle$
 where $m^{\tau}_{t}$ is the number of $F\in \tau$ with arity $t;$
 and we call $\bar m^{\tau}$  the signature of $\tau$, of course when
 saying   ``the signature of $M$'' we mean ``of $\tau(M)
$''.

 (3) For $M$  a fim we call $B \subseteq   M$ closed in $M$
 (or $M-$closed )  if for
$b= F^{M} (a_1,\dots,a_{s})$  we have
 $b \in B$   iff   $a_1,\dots,a_{s} \in M.$
Let the closure
 of $A$ in $B$ or  cl$_{M} (A)$ for $A \subseteq M$,  be the minimal $M-$closed
set $B \subseteq$
 $M$ which include $A$. A close (non empty) subset
of $M$ is actually a submodel.
 We do not strictly distinguish
 between a closed subset $B$ of $M$
 and the model $M \restriction B$
(which are fims with the same vocabulary).

 (4) For $\tau-$index models $M,N$ let
 PHom$(M,N)$ be the set of
 functions $f$ from $P^{M}$ into
 $P^{N}$ such that
$x \le^{M} y \Leftarrow f(x) \le^{N} f(y).$
 Let Hom$(M,N)$ be the set of functions $f$ from $M$ into $N$ such that
 $f \restriction P^{M}  \in$ PHom$(M,N)$  and $b=F^{M}(a_1,\dots,a_{t})$
implies
 $f(b) = F^{N}(f(a_1),\dots,f(a_{t})).$
Let PHm$(M,N)$ be the set of functions
$f$ from $P^M$ into $P^N$, and let
Hm$(M,N)$ be the set of functions
$f$ from $M$ into $N$ \st
$f \rt P^M \in PHm(M,N)$ and
$b= F^M (a_1,\dots , a_t )$ implies
$f(b) = F^N ( f(a_1 ) , \dots , f(a_t ) )$.

 (5) Let Sort$^{M}(F)$  be the range of $F^{M}.$

 \

 1.2 Fact: (1) For any $f\in$ PHom$(M,N)$
 there is a unique $\hat f \in$ Hom$(M,N)$
 which extend $f.$

(2) For any $f\in$ PHm$(M,N)$
 there is a unique $\hat f \in$ Hm$(M,N)$
 which extend $f.$

 \

 1.3  CLAIM/DEFINITION:
 (1) For any fim $M$ there is a
polynomial {\bold
p}$(x)$, with
 rational coefficients but positive integers as values for
 $x$ a positive integer (really sum of binomial coefficients
 binom$(x,\lk m_1,...,m_n \rk)$  for $m_i =1,\dots$, arity$(\tau) )$
such that
 for $u \subseteq P^{M}$, the set  cl$_{M}(u)$ has  exactly {\bold p}$(|u|
)$ members.
 Now {\bold p}$(x)$ depend on the signature of $\tau$ only   and so we
shall denote it by
 {\bold p}$_{\tau} (x)$  or {\bold p}$_{M} (x).$ Note that {\bold
p}$_{\tau} (0)  = 0.$

 \

 1.4 DEFINITION:
 (1) We say that $\tau$ is  canonical vocabulary
 for $t$ (or $t-$canonical)  and write $\tau  = \tau_{t}$
 if $\tau = \{ F_1,\dots,F_{t},P, <\}$
 where arity$(F_{s})$  is $s.$

  (2) We say that $M$ is a $(J,t)-$canonical  fim  if:

$(a)$ $J$ is a finite linear order

 $(b)$  $M$ is a fim with the $t-$canonical vocabulary

 $(c)$  $(P^{M}, <^{M} )$ is $J$

 $(d)$  $F^{M}_1$  is the identity on $P^{M}$

 $(e)$ for $r=2,\dots,t$  the function $F^{M}_{r}$  is
 $F^{M}_{r} (a_1,\dots,a_{r}) =^{df}  \{ a_1,\dots,a_{r}\}.$

 \

 1.5 DEFINITION:
 (1)  Let  $M$ be a fim  with vocabulary
 $\tau =\tau_{M}$
 and let $\{ A_1, A_2\}$ be a partition of $P^{M}$
 to  convex sets  such that $A_1 <^{M} A_2$
which means that
$ (\pa a_1 \in A_1 ) ( \pa a_2 \in A_2 ) [ a_1 <^M
a_2] $.
  We define a
 vocabulary $\tau_{M,A_1,A_2}.$
 It contains, in addition
 to  the symbols  $P$,
 $<$,  for each function
 symbol $F$ of $\tau$ and a  $\le^{M}-$increasing sequence $\bar a_1$
  from $A_1$  and   a  $\le^{M}-$increasing sequence $\bar a_2$
 from $A_2$
 such that lg$(\bar a_1) +$lg$(\bar a_2) <$ arity$^{\tau}(F)$
 a function symbol called
  $F_{  \bar a_1, \bar a_2   }$ with arity
 arity$^{\tau}(F) -$ lg$(\bar a_1 ) -$ lg$(\bar
a_2).$

 We identify $F \in \tau$ with
$F_{ \langle \rangle, \langle \rangle }$
 and so consider $\tau_{M,\bar a_1, \bar a_2}$
 an extension of $\tau.$

 (2)  Let $\bar m =  \bar m [{\bold{ \tau}} ,k_0, k_1 ]$
be $\bar m  [{\bold{ \tau}}_{M,A_0,A_1}]$, the signature of
$\tau_{M,A_0,A_2}$ whenever $M$ is a $\tau-$fim of dimension $k_0 + k_1$
and $A_0$ is the
 set of $k_0$ first members of $P^{M}$ and $A_1$ is the set of $k_1$
 last members of $P^{M}.$

 (3) Let $M^{\tau}_{k}$  be a fim of vocabulary $\tau$
 and dimension $k$, say $P^{M} = \{1,\dots,k\}.$
 Let $\tau^{[k,\ell ]}$ be $\tau_{M^{\tau}_{k+\ell}, A_0,A_1}$
 where $A_0$ is the set of the first $k$ members of $P^{M}$
 and $A_1$ is the set of the last $\ell$ members of $P^{M}.$

 \

 1.6 DEFINITION:
 (1) Let {$\bar \Lambda$} denote a sequence $\langle \Lambda_{F} : F \in
\tau \rangle$
 where $\Lambda_{F}$ is  a finite alphabeth, and we let $\tau[{\bar
\Lambda} ] = \tau$, as {$\bar \Lambda$}
 determine $\tau.$  We call {$\bar \Lambda$} an alphabet sequence  (for
$\tau)$
 or a $\tau-$alphabet sequence.
 We may write $(\tau,\Lambda)$  instead {$\bar \Lambda$} if
 $\tau  = \tau[{\bar \Lambda} ]$  and $\Lambda_{F} = \Lambda$ for every $F
\in \tau.$

 (2)  We say $p$ is a  {$\bar \Lambda$}-type   if  $p$ is a function with
domain $\tau$
 such that $p(F)  \in \Lambda_{F}$;
 let {\goth p}, {\goth q} denote non empty sets of
 {$\bar \Lambda$}-types; we identify them with their
 characteristic functions that they define, so
  we assume that from {\goth p} we can reconstruct {$\bar \Lambda$}
 hence $\tau [ {\bar \Lambda} ].$
 Let {\goth p}$_{\bar \Lambda}$ be the set
 of {$\bar \Lambda$}-types. We may
write $\Lambda$  instead  $\bar  \Lambda$
 is $\Lambda_{F} = \Lambda$ for every
 $F \in \tau$ and then   let {\goth
p}$_{\tau,\Lambda}$ be the set of
 constant $(\tau,\Lambda)-$types.

 \

 1.7  DEFINITION:
 (1) For {$\bar \Lambda$}  a $\tau-$alphabet sequence,  let  $V=$
Space$_{\bar \Lambda} (M)$
 be defined as follows:

 its set of elements is  the set of functions $\eta$ with domain $M$,
 such that   $b \in$ Sort$^{M}  (F) \Rightarrow \eta ( b ) \in  \Lambda_{F} ;$
  we assume that from $V$ we can reconstruct $M$ and
$\Lambda.$

 (2) We say  $d$ is a  $C-$colouring of  $V$, if $d$ is a function form $V$
 into $C$, we say $c-$colouring if $C$ has $c$ members
and the default value
 of $C$ is [0,$c) = \{ 0,1,\dots, c-1\}.$

 (3) We say $L$ is a $V-$line  or a line of $V$  if for {\goth q} $= {\goth
p}_{\bar \Lambda}$
 we have :    $L$ is a $(V, {\goth q} )-$line  or a {\goth q} -line of $V ;$
 this is defined
 for
 {\goth q}   a  (non empty) subset of {\goth p}$_{\bar \Lambda}$
 and it  means:

  $L$ is a subset of $V$ such that for some subset
  supp$(L) =$ supp$_{M} (L)$
 of $M$ we have:

 (a) supp$(L) \cap P^{M}$ is non empty  and we call it supp$^{P} (L)$

 (b)  supp$(L)$ is the $M-$th closure of supp$^{P}(L)$

 (c) for any $\eta,\nu \in L$ we have
 $\eta \restriction (M  \setminus$ supp$_{M}(L))  = \nu \restriction (M
\setminus$ supp$_{M}(L) )$

 (d) for any $\eta \in L$  for some  $p \in {\goth q}$ we
 have  : if  $b \in$ supp$_{M} ((L)$ then  $\eta (b) = p( F_{M,b} )$

 (e) For any $p \in {\goth q}$ there  is $\eta \in L$
as in clause $(d).$

 (5)  For $L$ as above
 and $p  \in  {\goth p}$ let pt$_{L} (p)$ be the unique $\nu \in L$
 such that for every  $a$ $\in$ supp$_{M} (L)$ we have $\nu(a)= p( F_{M,b}).$
 For {\goth q}$^* \subseteq {\goth q}$, the {\goth q}$^* -$ subline of a
{\goth q}-line $L$
 is $\{$pt$_{L} (p) : p \in {\goth q}^* \}.$

 (6) For a colouring $d$ of $V$, we say a $V-$line
  (or $(V,{\goth q} )-$line)  $L$ is $d-$monochromatic
 if $d$ is constant on $L.$

 (7) When  we are given $M,\tau,{\bar \Lambda},V$ as in
 part (4)
 and in addition we are given $m$, we define when $S$ is an
 $m-$dimensional  $V-$subspace , or $m-$dimensional
 subspace for $V.$
  It means that for some  sequence $\langle M_{\ell} :\ell<m\rangle$
 we have

 (a) each $M_{\ell}$ is a  submodel of $M$,

 (b) if  $\ell_1 < \ell_2 <m$
  then $M_{\ell_1}, M_{\ell_2}$ are disjoint,

 (c) for some $\rho$,  a function with domain
  ($M\setminus$  cl$(\cup\{M_{\ell}:\ell<m)\}))$
   such that $\rho(b) \in \Lambda_{F_{M,b}}$
for every $b \in$ Dom$(\rho)$,
 and some $m-$dimensional $\tau-$fim $K$
say $K = M^{\tau}_{[0,m)}$
  and letting $N$ be the submodel of $M$
 with universe  cl$_{M}( \bigcup_{\ell < m} M_{\ell})$
{\bold  there is }
  $f \in$ Hm$( N, K)$   which is onto $K$
 such that $f \restriction P^{M_{\ell}}$ is
constant for each $\ell$  and: $\nu \in S$ iff
$\nu$ extend $\rho$
and  for some $\varrho \in$ Space$_{\bar \Lambda}(K)$
 we have $b \in N \Rightarrow \nu(b) = \varrho (\hat f (b))$

 (8) We call $S$ convex if

(a) for $\ell_1<\ell_2  <m$
and $a_1 \in M_{\ell_1}$
 and $a_2 \in M_{\ell_2}$  we have
$a_1 <^{M}  a_2$
and

(b) $ f \in $ Hom$(N,K)$.

 (9) For $S$ as above and (see Definition 1.5(3))
 $\varrho \in$ Space$_{\Lambda} ( M^{\tau}_{[0,m)})$ we define
 pt$_{S} (\varrho)$ as
  the unique $\nu\in S$ as above in part (8).

 \

 We may define now  a natural function,
 which is our main concern here :

 1.8 DEFINITION: (1) Let $f^1({\goth p},c)$ where {\goth p} $\subseteq$
 {\goth p}$_{\bar \Lambda}$   ( and {$\bar \Lambda$} an alphabet sequence)
be the minimal $k$ such that for any
 $\tau_{\bar \Lambda}-$fim $M$ of dimension $k$,   we have:

 for any $c-$colouring $d$ of $V =$ Space$_{\bar \Lambda} (M)$ there is a
{\goth p}-line $L$ of $V$
  which is $d-$monochromatic,  i.e. such that $p, q \in {\goth p}$ implies
 that pt$_{L} (p )$ , pt$_{L} (q)$  have the same colour
 (by $d).$

  If $k$ does not exist we may say it is $\omega$
 or is $\infty.$ We may write
 $f^1_{\tau}({\goth p},c)$ or $f^1({\goth p},c;
\tau)$  to stress the role of $\t$.

 (2) If {\goth p} $= {\goth p}_{\bar \Lambda}$
  we may write $f^1({\bar \Lambda},c).$
   If  $\Lambda_{F} = \Lambda$ for every  $F \in \tau = \tau [{\bar
\Lambda} ]$
 then we may write   $f^1_{\tau}(\Lambda, c) ;$
 in this case we can replace $\Lambda$ by $|\Lambda|.$
 Clearly only $\bar m^{\tau}$ is important
 so we may write only it. Also we may write
 $f^1_{\tau} ( \bar n, c)$ for $f^1_{\tau}({\bar \Lambda}, c)$ whenever
  $\bar n = \langle n_{F} : F \in \tau \rangle$
 and $|\Lambda_{F}| = n_{F}.$

 \

 We can of course  use the multidimensional
 versions of those definitions

 1.9 DEFINITION: Let $f^1({\bar \Lambda},m,c)$ be the minimal $k$ such that
for any
 $\tau-$fim $M$ of dimension $k$  we have:
 for any $c-$colouring $d$ of
Space$_{\bar \Lambda}(M)$ there is a
 convex  subspace $S$ of $V$ of dimension $m$
  which is $d-$monochromatic,  i.e. such that all the points
 in $S$ have the same colour (by $d)$,
 if $k$ does not
  exist we say it is $\omega$
 or is $\infty.$  We may write
$f^1_{\tau} ({\bar \Lambda}, m,c)$
 or $f^1_{\tau} (\Lambda,m,c)$  etc. as before.
  Clearly only $\bar m^{\tau}$ is important
(rather than $\t$),
 so we may write only it.
 We may replace $\Lambda$ by $|\Lambda|.$
 We may replace {$\bar \Lambda$ } by
 $\langle n_{F} : F \in \tau \rangle$  when
 $n_{F} = |\Lambda_{F} |.$

 At present, it does not really matter
  if we omit the demand "convex¾ above.

 \

 The function has some  obvious monotonicity properties,
 we mention
 those we  shall actually use.

 \

 1.10 Claim:
 (1) For $\ell= 1,2$  assume
  {$\bar \Lambda$}$^{\ell}$ is an alphabet sequence for the
 vocabulary $\tau^{\ell}$ and arity$(\tau^1 ) \le$ arity$(\tau^2)$
  and for each $m= 1,\dots$, arity$(\tau^1)$  we have

 $\Pi\{|\Lambda^1_{F}| : F \in \tau^1$ has arity $m \} \le$
 $\Pi\{|\Lambda^2_{F}| : F \in \tau^2$ has arity $m \}.$

 {\bold Then}  $f^1({\bar \Lambda}^1,c) \le f^1({\bar \Lambda}^2, c ).$

 (2)  For $\ell= 1,2$  assume
  {$\bar \Lambda$}$^{\ell}$ is an alphabet sequence for the
 vocabulary $\tau^{\ell}$ and $\tau^1 \subseteq \tau^2$
 and {$\bar \Lambda$}$^1 = {\bar \Lambda}^2 \restriction \tau^1$
 and $F \in \tau^2 \setminus
 \tau^1 \Rightarrow |\Lambda^2_{F}| = 1.$

 {\bold Then}  $f^1({\bar \Lambda}^1,c) =
 f^1({\bar \Lambda}^2, c ).$

Proof: Straightforward.

$\qed_{1.10}$

1.11 DEFINITION: We define, for $\ell= 1,2,3$ what is a
fim$^{\ell}$,  we just replace in Def.$1.1$
 clauses $(d),(e)$ by

 $(d)_{\ell}$ $F^{M}$ is a partial function  such that if
$F^{M}(a_1,\dots,a_{r})$
 is well defined  (so $r =$ arity$(F) )$ then
  $a_1,\dots,a_{m} \in P^{M}$
 and $\ell= 1$ implies the function is symmetric, i.e. does not depend
 on the order of the variables, so if not said otherwise we assume
 $a_1 \le^{M}  a_2  \le^{M} \dots \le^{M} a_{r} .$

 $(e)_{\ell}$  if $F^{M}_1(a_1,\dots,a_{r}) = F^{M}_2 (b_1,\dots,b_{t})$
 and $\ell \in \{1,2\}$
 then $F_1 = F_2$ (hence $r=t )$  and
 $\ell = 2 \Rightarrow \bigwedge_{s=1,\dots,r} a_{s} = b_{s}$

 and

  $\ell = 1 \wedge \bigwedge_{s=1,\dots,r-1} a_{s} \le^{M} a_{s+1}$
 $\wedge \bigwedge_{s=1,\dots,r-1} b_{s} \le^{M} b_{s+1}$
 $\Rightarrow \bigwedge_{s=1,\dots,r}a_{s} = b_{s}.$
  So we  let
  base$_{M} (b) =^{df}   \{ a_1,\dots,a_{r}\}$
  and  when $\ell=1,2$  let
base$_{s}(b) =$ base$_{M,s}(b) =^{df} a_{s}$
 where $b = F^{M} (a_1,\dots,a_{r})$
 (and  if $\ell = 1$  then  $a_1 \le^{M}  a_2  \le^{M} \dots \le^{M} a_{r}$
, of course)
 and   $F_{M,b}  =^{df}  F ;$ those are well defined by the
 demand above.

\noindent
 $(e)_{\ell}^{'}$  if ${\ell} \in \{ 1,2,3 \} $ and $b \in M \setminus  P^{M}$
 then for some $F\in \tau$ and $a_1,...,a_{arity[F]} \in P^M$ we have $b =
F^M(a_1,...,a_{arity[F]})$.
 So $\ell =1$ is the old notion and for $\ell =3$ we require very little.
We define $f^{\ell}_{\lambda}({\bar \Lambda}, c)$ as in
 Definition 1.9 for fim$^{\ell}$ (so again $\ell =1$ is our standard case.)

\

  1.12 Claim: Let $\tau$ be a vocabulary and $\tau_{\circ} = \{
G_{F,\pi} :$
 $F \in \tau$ and $\pi$ is a permutation of $\{1,\dots$,arity$(F)\} \}$

 with arity$(G_{F,\pi}) =$ arity$(F).$

 {\bold Then}

 $(\alpha)$   If {$\bar \Lambda$} is a $\tau-$alphabet sequence and
  {$\bar \Lambda$}$^{\circ} = \langle \Lambda^{\circ}_{G} : G\in
\tau_{\circ} \rangle$
 where $\Lambda^{\circ}_{G_{F,\pi}}  = \Lambda_{F}$
 then  {{$f^2$}}$_{\tau} ({\bar \Lambda},c) \le f^1_{\tau_{\circ}}({\bar
\Lambda}^{\circ},c)$

 $(\beta)$ For {$\bar \Lambda$} a
$\tau-$alphabet sequence we have:
 $f^3_{\tau}({\bar
\Lambda},c)$ is at most
 RAM$\(f^2_{\tau}({\bar \Lambda}^{\circ},c \)$,
arity$(\tau), c^*)$
 where  e.g. $c^*$ depend on $\tau$ only (and RAM
stand for Ramsey number).

 $(\gamma)  f^1_{\tau}({\bar \Lambda},c) \le f^2_{\tau}({\bar \Lambda},c)$

 $(\delta) f^2_{\tau}({\bar \Lambda},c)\le f^3_{\tau}({\bar \Lambda},c)$

 \

 Proof:
 Straightforward.

 $\qed_{1.12}$

 \page
 {\bold Section 2 : Proof of the partition
  Theorem with a bound}

Except Def 2.1,2.2 this section is for the
reader convenience  only, as it give
a proof of a weaker version of
the first theorem  (with a bound which we get
 by triple induction).
Later
in 4.1-4.10
we give a complete proof with the
primitive recursive bound, formally not depending
on the proofs here. The strategy is to make the $b \in M$
with $|{\rm base}_M(b)|$ maximal immaterial. We first define some
help functions.

 2.1 DEFINITION:
 (1) We call a vocabulary $\tau$ monic
 if there is a unique  function symbol
of maximal arity,
 we then denote it by $F^{{\rm max}}_{\tau}.$

 (2) For a $\in P^{M}$ let $M_{a}$ be cl$_{M} (\{ P^{M} \setminus  \{a\} )$

 (3) For $V =$ Space$_{\bar \Lambda} (M)$ and  $N$ a
closed subset of $M$
 and $H \in \tau$,
 we say that a colouring $d$ of $V$ is $(N,\alpha,H)-$invariant
 if : $\alpha \in \Lambda_{H}$, and the following
holds, for any  $a \in P^{N} :$

 (*) if $\nu,\eta \in V$  and
$\nu \pr  M_{a}  = \eta \pr M_{a}$
and
 $\[ b \in M \wedge
 {\rm base}(b) = \{a\} \wedge F_{M,b} = H  \Rightarrow
\nu(b)= \alpha =  \eta(b) \]$
 then $d(\nu) = d(\eta).$

 (4) In part (3) we write
$(\ell,\alpha,H)-$monochromatic if above
 $N$ is such that $P^{N}$ is the set
of the last $\ell$ members of $P^{M}.$
We write $(M, \a , H)-$monochromatic if in part (3)
we have $M=N$.

 (5) In parts  (3) and (4) we may omit $H$
 when $\tau$ is monic and
 $H = F^{{\rm max}}_{\tau}.$ Replacing $\alpha$
 by $\Lambda^*$ mean that $\Lambda^*$
 is a subset of $\Lambda_{H}$ and the demand
holds for every $\alpha \in \Lambda^*.$

 \

 2.2 DEFINITION:
  Let $f^0$ be defined as follows.
First, $f^0_{{\bar \Lambda}} (n, \ell,
c)  =$
 $f^0_{\tau,{\bar \Lambda}} (n, \ell, c)$  is defined iff
 $\tau$ is monic with $H =  F^{{\rm max}}_{\tau}$
and {$\bar \Lambda$}   an alphabet sequence
 for $\tau$
 and $n \le |\Lambda_{H}|$  and $n < |\Lambda_{H}| \vee
 ( n=|\Lambda_{H}| \wedge \ell =0)$.
 Second, $f^0_{{\bar \Lambda}} (n, \ell,
c) $ is
 the first $k$ (natural number, if not
defined we can understand it as $\infty$ or $\omega$ or `does not
exist` )  such that (*)$_{k}$  below
 holds,

  where:

 (*)$_{k}$ If clauses (a)-(f) below hold then
 there is a $d-$monochromatic line  of $V$, where :

 (a)  $M$ is a fim of vocabulary $\tau$

 (b) the dimension of $M$ is $k$

 (c)  $V =$ Space$_{\bar \Lambda} (M)$

 (d) $\Lambda^{\circ}$ is a subset of $\Lambda_{H}$
with exactly $n$ members

 (e) $d$ is an $(M,\Lambda^{\circ}, H)-$invariant
colouring of $V$

 (f) if $\ell \not= 0$, then
there is an $\a$ such that
 $\alpha \in \Lambda_{H}
\setminus  \Lambda^{\circ}$
 and $d$ is $(\ell, \a , H)-$invariant.


 \

 Immediate connections are:

 2.3 Observation:
 (1) The function
$f^0_{\tau,{\bar \Lambda}}(n,\ell,c)$ increases
 with $c$
 and decreases with $\ell$  and $n.$

 (2)
 The function $f^0_{\tau ,{\bar \Lambda}}
(n,\ell ,c)$
  depends just on $n,\ell ,c$ and the
 set $\{$(arity$(F),|\Lambda_{F}|) :F \in \tau\}$
 (possibly with  multiple membership),
  so we may replace $\tau$ by its ${\bar m}^\t $
 (similarly for other such functions).

 (3) In definitions 1.8,1.9,2.2 the demand holds for any larger $k.$

 (4) $f^0_{\bar \Lambda} (0,0,c) = f^1({\bar \Lambda},c).$

 (5)  If $\tau$ is monic and
$H = F^{{\rm max}}_{\tau}$
 and $\tau^- = \tau \setminus \{H\}$  then
 $f^0_{\bar \Lambda}(|\Lambda_{H}| , 0 ,c) =
 f^1({\bar \Lambda} \restriction \tau^-, c).$

 (6) If $\ell^* = f^0_{\bar \Lambda}(n+1,0,c)$ then $f^0_{\bar \Lambda}
(n,\ell^*, c) = \ell^*.$

 Proof: Trivial.

 \

 2.4  MAIN  Claim:
   Assume

 (a) {$\bar \Lambda$} is an alphabet sequence
for  a vocabulary $\tau = \tau [{\bar \Lambda} ]$,
 and $n < |\Lambda_{H}|$

 (b) $\tau$ is a monic vocabulary with
$H = F^{{\rm max}}_{\tau}$

 (c)  $k_0 \ge  f^0_{{\bar \Lambda}}(n, \ell +1, c)$
  and $k_0 > \ell$

 (d)  $K$ is a $\tau-$fim  of dimension $k_0  - 1$
 with $A_2$ the last $\ell$ elements and
$A_1$ the first
 $(k_0 -\ell -1) -$elements  (this $K$ serve
just for notation)

 (e)  $\tau^*$ is the vocabulary
  $(\tau_{K,A_1,A_2}) \setminus \{H\}$
 see Definition 1.5(3);
  so

(i) arity$(\tau^*) <$ arity$(\tau) )$,

(ii) proj is the following function
 from $\tau^*$ to $\tau:$ it map
$F_{K,\bar a_1, \bar a_2}$
 to  $F$
so proj$\restriction \tau$ is the identity, and

(iii) ${\bar \Lambda}^*
=^{df} \langle \Lambda^*_{F} : F \in \tau^*
\rangle$
 where $\Lambda^*_{F} = \Lambda_{{\rm proj}(F)}$.

  (f) $c^* =^{df} c^{{\rm card}({\rm Space}_{\bar
\Lambda}(K))}.$

  {\bold Then}

 $f^0_{\bar \Lambda}(n,\ell,c)  \le k_0 + $
  $f^1({\bar \Lambda}^*, c^*  )  -1$

 \

 Proof:
 Let  $k_1 = f^1 ({\bar \Lambda}^*, c^*  )$
 and let $k= k_0 + k_1 - 1$, so it suffice to prove
 that $k \ge  f^0_{\bar \Lambda} (n ,\ell , c ).$
 For this it is enough to check (*)$_{k}$
from Definition 2.1(1),
 so let  $\Lambda^\circ$  be a subset of
 $\Lambda_{H}$ with $n$
 elements and
$\alpha^* \in \Lambda_{H} \setminus \Lambda^*$,
also  let
 $M$  be a fim of vocabulary $\tau$
 and dimension  $k$ (i.e. $P^{M}$ is  with  $k$
members), $V =$ Space$_{\bar \Lambda} (M)$,
 and $d$ an  $(\ell, \alpha^*,H)-$invariant
 and $(M,\Lambda^\circ,H)-$invariant  $C-$colouring   of $V$
 such that $C$ has  $\le  c$ members.  So we just have to prove  that
  the conclusion of Definition 2.2 holds, which means
 there is a $d-$monochromatic line of $V$.

 Let
 $w_1 =^{df}  \{ a: a\in P^{M}$
 and the number of $b<^{M}$ $a$ is $\ge k_0  - \ell - 1$ but
  is $<k_0 - \ell - 1 +k_1\}$
 hence in $w_1$  there are $k_1$ members,  and let $w_0$ be the set of
 first  $k_0  -\ell -1$  members of $P^{M}$ by  $<^{M}$,
and  lastly
 let $w_2$ be the set of
 the $\ell$ last members of $M$ by $<^{M}.$  So $w_0, w_1, w_2$  form a
convex
 partition of $P^{M}.$

 Now  we let
  $K$ be $M$
 restricted to cl$_{M} (w_0\cup w_2)$,
 (note that this gives no contradiction to the
assumption on
$K$ i.e. clause (d) of the assumptions,
 as concerning  $K$ there,
only its
 vocabulary and dimension are important and they fit).
  Let   $K^+$ be a fim
 with vocabulary $\tau$  and dimension $k_0$,
 let $g_0 \in$ PHom$(M,K^+)$ be the
 following function from $P^{M}$ onto $P^{K^+}:$
 it maps all the members of $w_1$ to one member of $P^{K^+}$
 which we call $b^*$, it is a one to one order preserving function
 from $w_2$  onto
 $\{ b \in P^{K^+} : b^* <^{K^+} b \}$
 and it is a one to one order preserving function
 from $w_0$  onto
 $\{ b \in P^{K^+} : b <^{K^+} b^* \}.$
 Let $g  \in$ Hom$(M,K^+)$ be the unique
 extension  of $g_0;$
 without loss of generality $g_0$ is the identity on $w_0$ and on $w_2$
 hence without loss of generality $g$
is the identity on $K$, it exist by 1.2.

 Next recall that the vocabulary
$\tau^* = \tau_{K, w_{o},w_2} \setminus \{H\}$
 is a well defined vocabulary ( see Definition
 1.5(1)
and remember that $\t \ps \t_{K,w_0 , w_2}$
so $H\in \t_{K , w_0 , w_2}$).  Next
 we shall define a $\tau^*-$model $N.$
 Its universe is  $M \setminus  K \setminus  A^*$
  where
 $A^* =^{df} \{ b \in M :$ base$_{M}(b) \subseteq w_1$ and
 $F_{M,b} = H \}$,  we let $P^{N}$  be $w_1$  and $<^{N}$  be  $<^{M}
\restriction P^{N}.$
 Now  we have to define each function $F^{N}_{K,\bar a_1, \bar
a_2}$,
 say of arity $r$, where $F \in   \tau, \bar a_1$ a non decreasing sequence
form $w_0$
 and $\bar a_2$ a non decreasing sequence from $w_2$,
 and lg$(\bar a_1)  +$
 lg$(\bar a_2 )  <$ arity(F)  and arity$(F_{K,\bar a_1, \bar a_2 } )$
 $<$ arity$(\tau).$  Note that the last condition
 is equivalent to : if $F=H$ then at least one of the sequences
 $\bar a_1, \bar a_2$ is not empty.

 For $b_1 \le^{N} \dots \le^{N} b_{r} \in P^{N}$
 we let  $F^{N}_{\bar a_1, \bar a_2}
 ( b_1, \dots, b_{t} )$   be equal to

 $b= F^{M}( \bar a_1, b_1,\dots, b_{t},
\bar  a_2)$
 $= F^{M} ( a^1_1, a^1_2,\dots, a^1_{{\rm lg}(\bar a_1 )}$,
 $b_1,\dots, b_{t}, a^2_1, \dots, a^2_{{\rm lg}(\bar a_2) } ).$

 It is easy to check that the number
of arguments is
right
 and also the sequence they form is
$\le^{M}-$increasing, so this is well defined and
belongs to $M$, but still we have to check that it
belongs to $N.$
 First note that it does not belong to $K$, as
if $b \in K$ then
  base$_{\lg(\bar a_1 )+1}(b) \in K$
 and it is just $b_1$ which belongs to
 $w_1$, contradiction.
 Second note that it does not belongs to $A^*$, this
holds
  as we have substructed $H$ when we have defined $\tau^*.$
 Lastly it is also trivial to note
 that every member of $N$ has this form.
 It is easy to check that $N$ is really a $\tau^*-$fim.

 We next let  $V^* =$ Space$_{{\bar \Lambda}^*} (N)$ and
let $C^* = \{g: g$ is a function from
Space$_{\bar \Lambda} (K)$ to $C\}$ and
 define a $C^*$-colouring $d^*$ of  $V^*.$
 For $\eta \in V^*$ let $d^*(\eta)$ be the following function
 from Space$_{\bar \Lambda} (K)$ to $C$, letting $\varrho$ be the
  function with domain $A^*$
 which is constantly $\alpha^*:$
 for $\nu\in K$ we let $\( d^*(\eta) \) (\nu) =
d(\eta \cup \nu \cup \varrho).$

 Clearly the function $d^*(\eta)$ is a
 $C^*-$colouring of
Space$_{\bar \Lambda} (K)$.
 How many such functions there are?
 The domain  has clearly   card(Space$_{\bar \Lambda}(K))$
 members,
  (we can get slightly less if $\ell >0$, but with no real influence).
 The range has at most $c$ members, so the number of such functions
 is at most
$c^{{\rm card}({\rm Space}_{\bar \Lambda}(K)) }$, a
number which we have called $c^*$.

  So $d^*$ is a  $c^*-$colouring.

 Now as we have chosen
$k_1 = f^1({\bar
\Lambda}^*,c^*)$
 we can apply Definition 2.2
   to $V^* =$ Space$_{{\bar \Lambda}^*}(N)$
and $d^*;$
 so we can find a $d^*-$monochromatic
 $V^*-$line and we call it  $L^*.$
  Let  $h$ be the  function
 from  $U =^{df}$  Space$_{\bar \Lambda} (K^+)$  to $V$ defined as follows:

 (*) $h (\rho)=\nu$ iff :

 (a) $\nu\in V, \rho\in U$,

 (b)  $\nu \restriction
 K =\rho \restriction   K$

 (c) if  $b \in  N \setminus$ supp$_{N}(L^* )$
 (see Def 1.7(3)) then $\nu(b) = \eta(b)$ for every
 $\eta \in L^*.$

 (d) if a $\in A^* \setminus$  cl$_{M}$( supp$_{N}(L^*))$ then
 $\rho(a) = \alpha^*.$

 (e) if a $\in$ supp$_{N}(L^*)$, (so a $\in  N$,
 $F_{N,a} = F_{K,\bar a_1,\bar a_2}$,
 base$_{N}(a) \subseteq$ supp$^{P}_{N}(L^*) )$,
  and $b   \in   K^+$,
 $F_{K^+,b} = F, b= F(\bar a_1, b^*,\dots,b^*,\bar a_2)$
 (with  the number of cases of $b^*$ being
 arity( $F_{K, \bar a_1,\bar a_2} ))$
  then
  $\rho(b)=\nu(a).$

 (f) if $a \in A^* \cap$cl$_{M}$(sup$_{N}(L^*))$
and $b \in$
  $K^+$ is $H(b^*,\dots,b^* )$
 then  $\rho(b) = \nu(a).$

 \

 Let the range of $h$ be called $S.$ Now clearly

 $\otimes_1$ $(\alpha)$ $h$ is a one to one function from  $U$  to $S
\subseteq V.$

 $(\beta)$ $S$ has  $|$Space$_{\bar \Lambda} (K^+ )|$ members

 $(\gamma)$  $S$ is a subspace of $V$
of dimension $k_0$, \st
$h(\rho) =$ pt$_{S}(\rho)$,
 see 1.7(7).

 \

 Now clearly

 $\otimes_2$ there is a $C-$colouring $d^{\circ}$ of $U$ such that:

 $d^{\circ} (\nu) = d(h(\nu))$  for $\nu \in U.$

 \

and

\

 $\otimes_3$ (a) $d^{\circ}$ is
$(K^+ ,\Lambda^*
)-$invariant

 (b) $ d^{\circ}$ is $(\ell +1, \alpha^*,
H))-$invariant

 [WHY? Reflect]

 \

 Applying the definition
 of $k_0  \ge f^0_{\tau , {\bar \Lambda}}
(n, \ell +1,c)$ ,
 that is Definition 2.2
 to {$\bar \Lambda$}, $\alpha^*,U,d^{\circ}$ we can
conclude that there is
 a $d^{\circ} -$monochromatic $U-$line $L^{\circ}.$ Let
 $L =^{df} \{ h(\rho) : \rho \in L^{\circ} \}.$
 It is easy to check that $L$ is as required.

 $\qed_{2.4}$

 \

As a warm up for the later bounds we prove:

 2.5  Theorem:
 (1)   The function $f^1_{\tau}({\bar \Lambda},c)$  is  well defined, i.e.
always get value, a natural number.

 Moreover has a bound which we have got by triple induction.

 (2) Similarly the function $f^0.$

 \

 Proof:
 (1)   The proof follows by induction,
 the main induction is on $t=$  arity$(\tau_{\bar \Lambda} ).$
 Now by observation 1.10(1)  without loss of generality $\tau$
 is monic, i.e. has a unique function symbol of arity $t$,
  called $H =^{df} F^{{\rm max}}_{\tau}.$  Fixing $t$, we prove by
induction on
 $s =  |\Lambda_{H}|.$

 CASE 0: $t= 1$

 This is Hales-Jewett theorem (on a bound see [Sh:329] and [GRS80])

 CASE 1: $t>1, s=1$

 By claim 1.10(2) we can decrease $t.$

 CASE 2: $t>1, s \ge 2$

 We note that $f^1({\bar \Lambda},c) = f^0_{\bar \Lambda}(0,0,c)$ by 2.3(4)
 so it is enough to bound the later one. But by 2.3(5) we know
 $f^0_{\bar \Lambda}(|\Lambda_{H}|,0,c) = f^1({\bar \Lambda}\restriction
\tau^-, c)$ where $\tau^- =^{df} \tau \setminus \{H\}$, but
 for the later one we have a bound by the
induction hypothesis
 on $t$ as arity$(\tau^- )  \le t$,  so we have a bound on
 $f^0_{\bar \Lambda}(|\Lambda_{H}|,0,c).$ By the last two sentences
together, it is
 enough to find a bound to $f^0_{\bar \Lambda}(n,0,c)$ by downward
induction on
 $n \le |\Lambda_{H}|$, and we have the starting
case : $n=  | \Lambda_{H}|$
 and the case $n=0$ gives the desired conclusion.
 So assume we know for $n+1$ and we shall do it for $n.$
 Let $\ell^* =^{df} f^0_{\bar \Lambda}(n+1,0,c)$,  so we know
 that $\ell^* = f^0_{\bar \Lambda}(n,\ell^*,c)$ by 2.3(6),
 so we by downward induction on $\ell \le \ell^*$
 give a bound to $f^0_{\bar \Lambda}(n,\ell,c).$
 So we are left with bounding
  $f^0_{\bar \Lambda}(n,\ell,c)$  given bound for
 $f^0_{\bar \Lambda}(n,\ell+1,c)$
  ( and also
 $f^1_{\tau_{\circ}} ({\bar
\Lambda}^{\circ},c^{\circ})$ whenever
arity$(\tau_{\circ} ) <t).$
 For this 2.4 was designed, it says

 $f^0_{\bar \Lambda}(n,\ell,c) \le
f^0_{\bar \Lambda}(n,\ell+1,c) + f^1_{\t^*}({\bar
\Lambda}^*,c) +1 $

 where $\tau^*,{\bar \Lambda}^*$ were
 defined there and arity$(\tau^*) <$ arity$(\tau)$,
(well, we have to assume that
$\ell < f^0_{\bar\L} (n , \ell +1 , c )$, but
otherwise  use $ \ell +1 +  f^1_{\t^*}({\bar
\Lambda}^*,c) +1 $

 $\qed_{2.5}$

 \page
{\bold  Section 3 : Higher Dimension Theorems}

\

Concerning the multidimensional case  (see Def
1.9):

 3.1 Conclusion:
 (1) For any {$\bar \Lambda$}, $m$ and $c$, we have
$f^1({\bar \Lambda},m,c)$ is well defined
 (with  bound as in the proof, actually
using one further induction  using only $f^1_{\t_i} (
{\bar \L} , c )$ for suitable $\t_i$-s in teh $i-$step.)

 (2) We can naturally
 defined $\tau-$fim of dimension
$\aleph_0$
 and convex subspaces, and prove that  for any
 $\tau-$fim $M$
 of dimension $\aleph_0$ and alphabet
sequence $\bar \Lambda$,
 if Space$_{M}(\bar \Lambda)$ is the
union of finitely many Borel subsets,
 then  some convex subspace  $S$ of dimension $\aleph_0$
 is included in one of those Borel subsets.

 Proof:
 (1) For simplicity (and without loss of
 generality by 1.10(1))  we have
 {$\bar \Lambda$} is constantly
 $\Lambda$, so each
 $\Lambda_F$ is  $\Lambda$,  a
fixed alphabet.
 We   choose by induction on $i = 0,\dots , m$
 the objects  $M_{i},  \tau_{i}, k_{i}$
and $c_i$ such that

 (a) $k_0 = 0$ and $k_{i} < k_{i+1}$

 (b) $M_{i}$  is a fim for $\tau$ of dimension
$k_{i}$ (we allow
 empty fim, if you do not like it start with $k_0 =1$)

 (c)  $M_{i+1}$ is an   end extension of $M_{i}$

 (d) $\tau_{i} = \tau_{M_{i}, P^{M_{i}}, \0 }$ (see
Definition 1.5(1) )

 (e) $c_0$ is $c$ and $c_{i+1}$ is
$c^{ |{\rm Space}_{\Lambda} ( k_{i} + m- i) | }$

 (f)   $k_{i+1} = k_{i}  +
f^1_{\tau_{i}}(\Lambda,c_{i}).$

 \

 There is no problem to carry the definition
 and we can  prove that $k_{m}  \ge f^1_{\tau}(\Lambda,m,c).$

 The proof is straight.

2) Such theorems are closed relatives to
theorems on appropriate forcing notions, as anyhow
it is a set theoretical theorem
we use such approach.
Specifically we use the general treatment of
creature forcing of [RoSh 470].
For any finite non empty $u \ps \w$ let
$M^\t_u$ be a $\t-$model with
$(P^{M^\t_u} , \le^{M^\t_u} )
= (u, \le)$, and without loss of
generality
 $u_1 \ps u_2 \A M^\t_{u_1} \ps M^\t_{u_2}$.
So for infinite
$u \ps \w$ we have
$M^\t_u = \pu \{ M^\t_{u_1}
: u_1 \ps u $ finite $\}$
is well defined.

A $\bar \Lambda$-creature
$\gc$ consist of
a convex subspace $S^{\gc} = S[{\frak c}]$
of some $M^\t_u$
for some finite non empty $u = u[ {\frak c}]$
of the form $[n,m] = [ n_{\frak c} , m_{\frak c} ]$
.

Fot creatures ${\frak c}_1, \dots , {\frak c}_k$ we let
$\S ({\frak c}_1, \dots , {\frak c}_k )$
be well defined
iff $m_{{\frak c}_\ell}  = n_{{\frak c}_{\ell+1}}$ for
$\ell \in [1,k)$
and it is the set
of ${\bar \Lambda}-$creatures ${\frak c}$
such that
$n_{\gc}  =  n_{\gc_1}, m_{\gc} = m_{\gc_k}$
 and $ \h \in S^{\frak c} \pc
\ell
\in [1,k)  \A
\h\rt u[{\frak c}_\ell ] \in S^{{\frak c}_\ell}$.

So the forcing notion $\xQ$
is well defined by
[RoSh 470]
for the case "the lim-sup of the norms is infinity''.
So a condition $p$ has the form
$\lk \h , \gc_1 , \gc_2 , \dots \rk =
\lk \h^p , \gc^p_1 , \gc^p_2 , \dots \rk$
where for $t= 1,2,\dots $, $\gc^p_t$
 is a $\bar \L-$creature , $m_{\gc^p_{t+1}} = n_{\gc^p_t}$.
Let ${\xB} =^{df}   {\rm Space}_{\bar\L} (M^{\t}_\w )
= \{ \r: \r \hbox{ is a function with domain }
 M^\t_\w \hbox{ satisfying } f(b) \in \L_{F(b)} \}$ where
$F(b) = F_{M^\t_{\w},b}$.
We say that $\r \in \xB$
obeys $p \in \xQ$ if $\h^p \ps  \r$ and for
$t=1,2,\dots $ we have $ \r \pr u^{\gc^p_t}  \in S^{\gc^p_t}$.
It is proved there that
such forcing notions has many good properties. In particular
letting
cont$(p) = \{ \r:  \r \in {\xB} \hbox{ obeys } p \}$
is a function with domain
and defining the $\xQ$-name
$  {\hat  f} = \pu \{ f^p : p \in {\hat G}_\xQ \}$.
Now note that

(a) $ p \pf_\xQ " {\hat f} \in$ cont$(p)$"

(b) if $N \pE (\sH (\chi ), \in )$
is countable, the definition
of those countably many Borel
sets belongs to $N$, and
$p \in \xQ \pv N$, then
we can find $q$ such that

(i) $p\le q$

(ii) every $f \in {\rm cont}(q)$
is a generic for $\xQ$
over $N$

(iii) for some $p' , n'$  we have
 $ p \le p' \in N \pv \xQ, p' \le q$
  and  $ p' \pf_\xQ "{\hat f} \in  A_{n'} "$

Together we conclude that
 ${\rm cont}(q) \ps A_{n'}$
and we are done.

 \

 $\qed_{3.1}$

 \

 We turn to relating the old results from Bergelson
Leibman [BL96]

 3.2 Conclusion:
 (1) Assume that

 (a) $\tau$  is a $t-$canonical vocabulary
 (see  1.4)

 (b) $k= f^1_{\tau} (\Lambda,c)$,  $\L$ a (finite) alphabet

 (c)  $R$ is a ring,  and $r_1,\dots,r_{k} \in R$

 (d) for $\alpha \in \Lambda, {\bold p}_{\alpha}(x)$ is a
 polynomial over $R$ (i.e. with parameters in
$R$).

 (e) $d$ is a $c-$colouring of $R$ ( actually enough
to consider
 a finite subset, the range of $g$
in the proof below)

 {\bold Then} we can find  $y,z$ and $w \subseteq
\{1,\dots,k\}$  such that

 $(\alpha)$ $y \in R$  and
 $z = \Sigma_{\ell \in w}
r_{\ell} \in R$

 $(\beta)$ the set
$\{ y + {\bold p}_{\alpha}(z) : \alpha \in \Lambda
\}$ is $d-$monochromatic

 (2) Assume   that

 (a) $\tau$  is a  vocabulary of arity $t$, such that for each $s = 1,\dots,t$
 in $\tau$ there are exactly $m^*$  function symbols of arity $s$

 (b) $k= f^1_{\tau} (\Lambda,c)$, $\L$ a (finite) alphabet

 (c)  $R$ is a ring,  and $r_1,\dots,r_{k} \in R$

 (d) for $\alpha \in \Lambda$ and
 $m<m^*,   {\bold p}_{\alpha,m}(x)$ is a
 polynomial over $R$ (i.e. with coefficients in
$R$).

 (e) $d$ is a $c-$colouring of
 $R^{m^*} = \{ \langle y_{m}:m<m^*\rangle  :
 y_0,\dots , y_{m^* -1} \in R\}$
  (actually enough to consider
 a finite subset, the range of $g$  in the proof below).

 {\bold Then} we can find  $y,z$  and $w \subseteq \{1,\dots,k\}$  such that

 $(\alpha)$ $y \in R$  and
$z = \Sigma_{\ell \in w}
r_{\ell} \in R$

 $(\beta)$ the set $\{ \langle y +
 {\bold p}_{\alpha,m}(z) : m<m^* \rangle :\alpha \in
\Lambda \}$ is $d-$monochromatic

 \

 Proof:
 (1)  Let $M$ be a fim for $\tau$ of dimension $k$
  and let $h$ be a  one to one order preserving
 function from $P^{M}$ to $\{1,\dots,k\}.$
 We define a function $g$ from $V=$ Space$_{\Lambda}(M)$  to $R.$
 For $\eta \in V$ we let $g (\eta) = \Sigma_{b \in M } g_{b}(\eta(b))$
 where
 $g_b$ is the following function from
$\L$ to $R$. For
 $b = F(b_1,\dots,b_{t}) \in M$  and $\alpha \in
\Lambda$
we let $g_b (\a)$ be
 zero if $\langle b_1,b_2,\dots,b_{t}\rangle$
is with repetitions
 and otherwise
 we consider
${\bold p}_{\alpha}(
\Sigma_{i=1,t} r_{h(b_{i})})$,
 expand it as sum of monoms
in $r_1,\dots , r_k$ , and let
 $g_{b}(\alpha)$ be the sum of
those monoms for
which
$\{ r_j : j \in \{ 1, \dots , k \} $ and $r_j$
appear  in
the monom $\}$ $=$ $\{ h( b_1 ) , \dots , h(b_t )\}$.
 Now we define a $c-$colouring  $d^*$ of $V$
by $d^* (\eta) = d(g(\eta)).$
 Let $L$ be a $d^*-$monochromatic line of $V$ ,
let supp$_{M}(L) = N.$
 Now let $y =^{df} \Sigma_{b \in  M \setminus N } g_{b}$(pt$_{L}(\alpha))$,
 note that all the $\alpha \in \Lambda$ gives the same value.
 Let $w =^{df} \{ h(b) : b \in$
supp$^{P}_{M} (L) \}$,
recalling Def 1.7(5)
 and so $z = \Sigma_{\ell \in w}r_{\ell}$,
now check.

\

 Note that algebraically it is more natural to
 defined $g$ differently,
 working by the rank of the monom rather
 that by the set of variables appearing.

 (2) Similarly, left to the reader.

 $\qed_{3.2}$

 \

 3.3 Discussion:          It is natural to ask:

 (1) Can we generalize the Graham Rothschild theorem?
(see [GR 71], [GRS 80])

 (2) Can we get here primitive recursive bounds?

 (3) Can we prove the density version of the theorem (2.11)?

 \

 Below we answer positively questions (1),(2),
 we believe that the answer to question (3)
 is positive too but probably it require
 methods of dynamical systems, see the book
 Furstenberg [Fu81].

 \

 3.4 DEFINITION:
 We define $f^4({\bar \Lambda},t,\ell,c)
= f^4_\t ( {\bar \L } , t, \ell , c )$
where 0 $\le
\ell < t$
 as follows. It is the minimal $k$ such that:
 if $M$ is  fim for $\tau, V =$ Space$_{\bar \Lambda}(M)$
  and $d$ is a $c-$colouring of  $\{ S: S$ is an $\ell-$subspace
 of $V \}$ {\bold then} for some
subspace $U$ of $V$ of dimension $t$,
 all the  $\ell-$subspaces of $U$
 (equivalently, $\ell-$subspaces of $V$
 which are contained in $U$)
have the same colour by $d.$

 \

 3.5 Theorem :
 (1) For any {$\bar \Lambda$}, $t,\ell, c$
as in Definition 3.3,  the function
 $f^4({\bar \Lambda},t,\ell,c)$ is well defined,
 i.e. is  finite.

 (2) Let $m =$ RAM$(t,\ell,c)$,
see  Definition
0.3(1),
where $\t$ is a vocabulary and ${\bar \L }$
is a $\t-$alphabet sequence,
  and define $k_{i}$ for $i=0,\dots,m$
 by induction on $i$ as follows (on $\tau^{[k,r ]}$
 see 1.5(3)):

  $k_0 = 0, {\bar \L}^0 = {\bar \L}$ and
$k_{i+1} = k_{i}  +
f^1_{\t_i} (\bar \Lambda^{i}, c_{i})$ where
$\t_i =^{df} \t^{[k_i , m-i ]}$ and
$\bar \Lambda^{i}$ is a $\tau^{[k,m-i ]}-$alphabet
sequence,  and
 $\Lambda^{i}_{F_{M^\t_{ k_i + m -i },  \bar a_1,
\bar a_2 }}$   has
$|
\Lambda^{i}_{F} | +  \ell + |M^{\tau}_{k_{i} +m-i}|$
 members and
$c_{i}   = c^{ {\rm card}({\rm
Space}_{\bar
\Lambda^{i} } (M^\t_{ k_i + m -i } ) )}.$

{\bold Then} $f^4_{\tau}(\Lambda,t,\ell,c) \le k_{m}.$

 Proof:(1) Follows from (2).

(2) Let $N = M^{\tau}_{\ell}$ (see notation in
1.5(3), recall that $\ell$ is the
 dimension of the subspaces we
are colouring)
 and let  $\{ \gamma_{a}:a\in N\}$ list
a set
 disjoint to $\Lambda$ without repetitions.

 We   choose  for  $i = 0,\dots , m$
 the objects  $k_i ,   \tau_{i} , {\bar
\L}^i$ (consistently with
 what is said in the statement of the theorem)
 and $M_i , M^+_i$,
 by induction on $i$ as follows:

\

  $\otimes_1$
 (a) $k_0 = 0$ and $k_{i} < k_{i+1}$

(b) $M_{i}$  is a fim for $\tau$ of dimension
$k_{i}$ (we allow
 empty fim, the space is the a singleton,
 if you do not like it start with 1)

 (c)  $M_{i+1}$ an   end extension of $M_{i}$
  and $M^+_{i}$ is an end extension of
 $M_{i}$ (so both have vocabulary $\tau)$
and has dimension $k_{i} +
m-i$

 (d) $\tau_{i} = \tau_{M^+_{i}, P^{M_{i}},
P^{M^+_{i}} \setminus P^{M_{i}} }$ (see Definition
1.5(3) )

 (e) ${\bar \L}^0 = {\bar
\L}$
and
$\Lambda^{i}_{F_{\bar a_1,
\bar a_2 }}$
  is the disjoint union of
$\Lambda_{F}, \Lambda^{\circ}_{F} =^{df} \{
\gamma_{b} : b \in N$ and
$F_{N, b} = F\}$     and
$\{ \beta_{b} :  b \in M^+_{i}$
 such that  $F_{M^+_{i}, b}  = F\}$
(and no two letter are incidentally equal,
of course).

 (f) $c_0$ is $c$ and $c_{i+1}$ is $c^{ {\rm
card } ({\rm Space}_{{\bar \Lambda}_{i}} (
M^\t_{k_{i}  + m- i}) )}$

 (g)   $k_{i+1} = k_{i}  +  f^1_{\tau_{i}}( \bar
\Lambda^{i},c_{i}).$

 \

 Let $k= k_{m}, M= M_{k}$  and let
 $V_{i} =$ Space$_{\bar \L}(M^\t_{i})$ and $V= V_{m}.$
 We shall regard an $\ell-$subspace $\Phi$  of $V$
as a function
 from $M$ to $\Lambda^{\circ} = \{ \gamma_{b}:b\in N\} \cup \Lambda$,  such
that (and where):

 $\otimes_2$
(a) $\Lambda = \bigcup_{F\in \tau}
\Lambda_{F}$,

 (b) $\Phi(b) \in \Lambda^{\circ}_{F_{M,b}} \cup
\Lambda_{F_{M,b}}$,  see clause (e) of $\otimes_1$

 (c) if $b \in M, \alpha \in \Lambda$  and
$(\forall \nu)[ \nu \in \Phi \Rightarrow \nu(b)=\alpha]$  then
$\Phi(b) = \alpha$

 (d) if $b \in M, a \in N$ and for every $\rho \in$
Space$_{\bar\Lambda}(N)$
 we have (pt$_{\Phi}(\rho))(b) = \rho(a)$ then
 $\Phi(b) = \gamma_{a}.$

 \

 (Reflect on the meaning of $\ell-$subspace of $M$, i.e. Definition 1.7(7)
 and it should be clear.)
 Let $d$ be a $c-$colouring of the set of $\ell-$subspaces of $V.$
 We shall define by downward induction on $i<m$ a  (non empty)
 subset $A_{i}$ of $P^{M_{i+1}}$  disjoint to $M_{i}$  and a function
$\varrho_{i}$
 from $B_{i} =^{df}  M \setminus$ cl$_{M} (M_{i} \cup
\bigcup_{j=i,\dots,m-1} A_{j}) \setminus  \bigcup_{j=i+1,\dots,m-1} B_{j}$
 into $\Lambda \cup \{ \beta_{a} : a  \in M_{i}^+\}.$
 We let  ${\bold R}_{i}$ denote the family of $\ell-$subspaces
 $\Phi$ of $V$ which satisfies:

 \

 (*)$_1$
 (a) if $j$ satisfies  $i \le j < m$  and
$b \in B_{j}$ and $\varrho_{j}(b) \in \Lambda$
 then $\Phi(b) = \varrho_{j}(b)$

 (b) if $j$ satisfies $i \le j <m$ and $b \in B_j$
and
$\varrho_{j}(b) = \beta_{a}$ where  a $\in M_{j}$
 then $\Phi (b) = \Phi(a)$

 (c) if $b_1,b_2$ satisfies the following then $\Phi(b_1) =\Phi(b_2)$
 where the demand is:

 (i) $b_1,b_2 \in$ cl$_{M} (M_{i} \cup
\bigcup_{j=i,\dots,m-1} A_{j})$
 and

(ii) $F_{M,b_1}  = F_{M,b_2}$ and   for every
  $r \in
\{1,\dots$, arity($F_{M,b_1} )$
  $\}$
 we have:
 base$_{M,r}(b_1) =$ base$_{M,r}(b_2)$ or they both
 belongs to the same   $A_{j}$
for some $ j \in
\{i,\dots,m-1\}$.

 \

 Now $A_{i},B_{i}, \varrho_{i}$ will be chosen such that the following
condition holds

 (*)$_2$ If $\Phi,\Psi  \in {\bold R}_{i}$ satisfy the clauses (a),(b) below
 then $d(\Phi)= d(\Psi)$  where

 (a) $\Phi \restriction$ cl$_{M}( M_{i} \cup
\bigcup_{j=i+1,\dots,m-1}A_{j}) =$
  $\Psi \restriction$ cl$_{M}( M_{i} \cup \bigcup_{j=i+1,\dots,m-1} A_{j} )$

 (b) if  $b \in N$ and $\gamma_{b} \in$ Rang$(\Phi \restriction M_{i+1})$
 then $\gamma_{b} \in$
 Rang$(\Phi \restriction M_{i})$.

 \

 Suppose now that we have carried this induction, and we shall show that this
 suffice. Let $S$      be the
 following subset  of $V$:

 (*)$_3$ $\eta \in S$ iff

 (a) if $i<m$ and $b \in B_{i}$ and
$\varrho_{i}(b) \in \Lambda$ then  $\eta(b) =
\varrho_{i}(b)$

 (b) if $i<m$ and $b \in B_{i}$   and
 $\varrho_{i}(b) = \beta_{a}$ and $ a\in M_{i}$ then
$\eta(b) = \eta(a)$.

 \

 Clearly $S$ is an $m-$subspace of $V$,
 and we may by (*)$_2$ above show that:

 (*)$_4$  if $\Phi$ is an $\ell-$subspace of $S$,
 then $d(\Phi)$ can be computed from
 $J[\Phi ] =^{df} \{$ Min $\{i: \Phi \restriction A_{i}$ is constantly the
$r-$th member
 of $P^{N}\} : r < \ell \}.$

 So for some function $e$, with domain
 the family of subsets
 of $\{0,\dots,m-1 \}$ with $\ell$ elements,
we have :
 if $\Phi$ is an $\ell-$subspace of
$S$ then  $d(\Phi ) = e(
J[\Phi] )$. Clearly the set ${\rm Rang }(e)$
has $\le \pl {\rm Rang (d) } \pl $ elements.

 By Ramsey theorem and the choice
of $m$, there is a subset $w$ of $\{0,\dots
, m-1\}$
 with $t$ members such that the
 function $e$ is constant on the family
 of subsets of $w$ with $\ell$ elements.
 Let $U$
be a subspace of $S$ of dimension
 $t$ such that if $b \in M$, base$(b)$
 not a subset of
 $\bigcup_{i \in w} A_{i}$ then
$\lk \n ( b ) : b \in U \rk$ is constant
(and the constant value belongs to $\L_{M,b}$.

  Clearly  $U$ is as required.
 The construction, i.e. the inductive choice
of $A_i , \R_i$ is straight.

 $\qed_{3.5}$

 \page

{\bold Section 4: The main Theorem}

 Now we turn to the obtainment of
 primitive recursive
bounds.
 The idea is that we decrease the dependency from
 below, dealing with the unary functions
 each time (rather than dealing with
 $H \in \tau$ of maximal
arity).

 \

 In the definition below, we shall use the case $r=1.$

 4.1 DEFINITION:
 (1) Recall that for $a$ $\in P^{M}$ we  let $M_{a}$ be cl$_{M} (\{ P^{M}
\setminus  \{a\} )$,
 that is $M$ restricted to this set.

 (2) For $V =$ Space$_{\bar \Lambda} (M)$ and  $N$ a
closed subset of $M$
 we say that a colouring $d$ of $V$ is $(N,r)-$base-invariant
 if the following holds, for any $a$ $\in P^{N} :$

 (*) if $\nu,\eta \in V$  and
$\n \pr M_a = \h \pr M_a$
 and $\[ b \in M \wedge r < |\{i:
i=1,\dots$, arity$(F_{M,b})$
 and  base$_{M,i}(b) = a\}|  \Rightarrow  \nu(b)=\eta(b) \]$
 then $d(\nu) = d(\eta).$

 (3) We write $(\ell,r)-$base-invariant if above
 $N$ is such that $P^{N}$ is the set of
the last $\ell$ members of $P^{M}.$

 \

 4.2 DEFINITION:
  Let $f^6$ be defined as follows. First,
 $f^6_{{\bar \Lambda}} (\ell, c)  =$
 $f^6 ({\bar \Lambda} , \ell, c) = f^6_{\tau}
({\bar \Lambda} , \ell, c)$  is defined iff
 {$\bar \Lambda$}   is an alphabet sequence
 for a vocabulary  $\tau$. Second, let
 $f^6_{{\bar \Lambda}} (\ell, c)  $  be
 the first $k$ (natural number, if not defined
 we can understand it as $\infty$ or $\omega$
or `does not exist`
)  such that (*)$_{k}$  below
 holds,
  where:

 (*)$_{k}$ If clauses (a)-(d) below hold then
 there is a $d-$monochromatic line  of $V$, where :

 (a)  $M$ is a fim of vocabulary $\tau$

 (b) the dimension of $M$ is $k$

 (c)  $V =$ Space$_{\bar \Lambda} (M)$

 (d)  $d$ is an  $(\ell, 1)-$base-invariant
colouring of $V$.

 \

 Immediate connections are:

 4.3 Observation:
 (1) The function   $f^6_{\bar \Lambda}(\ell,c)$ increases
 with $c$
 and decreases with $\ell.$

 (2) We have  $f^6_{\tau_1}( {\bar \Lambda^1}, \ell_1,c_1)$
 $\le  f^6_{\tau_2} (\bar \Lambda^2$, $\ell_2,c_2)$
 if:

(a) $c_1 \le c_2, \ell_1 \ge \ell_2$,   and

(b)  $s \le$ arity$(\tau_1)  \Rightarrow \Pi
\{|\Lambda^1_{F} | : F \in \tau_1$, arity$(F) = s\}
\le$
  $\Pi \{|\Lambda^2_{F} | : F \in \tau_2$, arity$(F) = s\}$  and

 (c) arity$(\tau_1) <s \le$ arity $(\tau_2) \wedge F
\in \tau_2 \wedge$ arity$(F)=s$
 $\Rightarrow |\Lambda^2_{F}| =1$

 (3) In definition  4.2  the demand holds for any larger $k.$

 (4) $f^6_{\bar \Lambda} (0,c) = f^1({\bar \Lambda},c).$

 Proof: Trivial.

 \

 4.4 Claim:
 Assume

 (a) $\tau$ is a vocabulary  of arity $>1$
  and $\bar \Lambda$ is a $\tau-$alphabet sequence

 $(b)$  $\tau^*$ is the following vocabulary: $\{ G_{F,e}: F \in \tau$,
   arity$(F) >1$ and $e$ is a convex equivalence relation
 on $\{1,\dots$, arity$(F)\}$ such that each $e-$equivalence class has
 at least two elements $\}$

 with arity$(G_{F,e}) =$ the number of $e-$equivalence classes
 and for some $H  \in \tau$
  of maximal arity, letting $e =^{df} \{(i,j): i,j \in [1$,arity$(H)]\}$
 we identify $G_{H,e}$ with id

 $(c)$  ${\bar \Lambda}^*$ is the following $\tau^*-$alphabet
sequence:
 $\Lambda^*_{G_{F,e}}  = \Lambda_{F} .$

 $(d)$   $\ell^* = f^1_{\t^*}({\bar \Lambda}^*,c)$.

 {\bold Then} $f^6_{\bar \Lambda} ( \ell^*,c) \le
\ell^*$.

 Proof:
 Let $M$ be a fim of vocabulary $\tau$ and dimension
 $\ell^*$ and $V=$ Space$_{\bar \Lambda} (M)$
 and $d$ is a $c-$colouring of $V$  which is $(\ell^*,1)-$base-invariant;
 it suffice to find  a monochromatic $V-$line $L..$

 Let $M^*$ be a fim of vocabulary $\tau^*$ and dimension $\ell^*$
 and $V^* =$ Space$_{{\bar \Lambda}^*}(M^*).$
 Let $g_0$ be an isomorphism from $(P^{M},<^{M})$ onto $(P^{M^*}, <^{M^*}).$
 We define a partial function $g$ from $M$ into $M^*$ as follows;
 if  $b=F^{M}(b_1,\dots,b_{t})$ so $t=$ arity$_{\tau}(F)$
 and  $b_1 \le^{M} b_2  \le^{M} \dots \le^{M} b_{t}$  and
 $e= \{(i,j): b_{i} = b_{j} \}$ and  $G_{F,e}$ $\in \tau^*$
  is well defined
 (i.e. every $e-$equivalence class has at least two elements)
 and the $e-$equivalence classes are
 [$s_{i},s_{i+1})$ for $i=1,\dots$,arity$(G_{F,e}) -1$ and 1$=s_1 <s_2<\dots$
 $<s_{{\rm arity}(G_{F,e})}=t+1$ then $g(b) =
G^{M^*}_{F,e}(g_0(b_{s_1}),\dots$,
 $g_0(b_{s_{{\rm arity}(G_{F,e})-1}})).$

 Note:

 $(*)_1$ $g$ is really a partial function from $M$ to
$M^*$

 $(*)_2$ if $\eta,\nu \in V$ and $\eta \restriction  {\rm Dom } (g) = \nu
\restriction {\rm Dom } (g)$ then
 $d(\eta) = d(\nu)$

 [Why? By the transitivity of equality,
it is enough to consider the case that for some
$a^* \in M \pb { \rm Dom}(g)$ we have
$\{ a^* \} = \{ a \in M : \h (a) \not=
\n (a) \}$.
Now by the definition of $g$
for some  $a \in P^M$
we have
$ ( \py ! i ) [ { \rm base}_{M,i} (b) =a ]$.
Now
 read
the Definition  4.1(2),(3) of the  "base invariant"]

 (*)$_3$  we can define a $c-$colouring $d^*$ of $V^*$
 such that $\eta \in V, \nu \in V^*$,  and
 $\[ b \in$ Dom$(g) \Rightarrow \eta(b) = \nu(g(b))\]$
 then $d(\eta) = d^*(\nu)$

 [why? by (*)$_2]$

 (*)$_4$ for any $V^*-$line $L^*$  there is a $V-$line $L$
 such that for every $\eta \in L$ for some $\nu \in L^*$ we have
 $d(\eta)=d^*(\nu)$

 [Why? Reflect. In details, let
 $w^* =$ supp$^{P} (L^*)$ and $N^* =$ supp$(L^*)$
 and $\nu^*$ is the function with domain
$M^* \setminus N^*$
 such that for every $b$ from this set and
$\nu \in L^*$ we have $\nu(b)= \nu^*(b).$
 Let $w =^{df} \{ b \in P^{M} : g_0 (b) \in w^* \}$
 and let $N =^{df}$ cl$_{M}(w)$ and
choose a function $\eta^*$ with domain
 $M \setminus  N$ such that for every
$b \in M \setminus N$ we have $\eta^*(b) =
\nu^*(g(b))$
 if $b \in{\rm Dom } (g)$  and is any member of
$\Lambda_{F_{M,b}}$  otherwise.
 Let $L$ be the $V-$line such that supp$(L)= N$ and for every $\eta \in L$
we have
 $\eta$ extend $\eta^*.$ Clearly $L$ is a
$V-$line and let $\eta \in L$ and we
  should check the desired conclusion.
So there is $p \in {\goth p}_{\bar \Lambda}$
 such that $\eta =$ pt$_{L}(p);$ now we define
$q \in  {\goth p}_{{\bar \Lambda}^*}$
 as follows: $q(G_{F,e}) = p(F)$, the later belongs to
 $\Lambda_{F}$  which is equal to $\Lambda^*_{G_{F,e}}.$ Let $\nu
=$pt$_{L^*}(q)$ and
 we should just check that
 $\eta,\nu$ are as in (*)$_3$ above so we are done.]

 By  the assumption $\ell^* = f^1({\bar \Lambda}^* , c))$
 (see clause $(d)$ in the assumption), hence there is a $d^*-$monochromatic
$V^*-$line $L^*.$
 Apply  (*)$_4$ to it, so there is a $d-$monochromatic $V-$line
 and so we are done.

 $\qed_{4.4}$

 \

 4.5 DEFINITION:
 (1) Assume the following:

(i) ${\bar \L}$ is an alphabet sequence for
 the vocabulary $\t$

(ii) ${\Bbb P}  \ps \{ ( p ,  q ):
{p} , { q}$ are ${\bar \L}$-types $\}$, see Def 1.6

(iii) $m,c> 0$.

We define $f^7_{\bar
\Lambda} ({\Bbb P},m,c)$ as the
  first  $k$ (if there is no such $k$
 it is $\omega$ or $\infty$ or undefined)  such that (*)$_{k}$ stated below
holds, where

 (*)$_{k}$ if clauses (a)-(e) below
 hold then there
is
 a subspace $S$ of $V$ of dimension $m$, satisfying:

 if $L$ is a $V-$line $\subseteq S$, and $(p,q ) \in {\Bbb P}$
  then $d$(pt$_{L}(p)) = d$(pt$_{L}(q))$

 where

 (a) $M$ is a fim of  vocabulary $\tau$

 (b) $M$ has  dimension $k$

 (c) $V =$ Space$_{\bar \Lambda} (M)$

 (d) ${\Bbb P}$ is a subset of $\{ (p,q): p,q \in
{\goth p}_{\bar \Lambda}$ and
 $\[ F \in \tau \wedge$ arity$(F)>1 \Rightarrow p(F) = q(F)\] \}$

 (e) $d$ is a $c-$colouring of $V$

 (2) Let {\Bbb P}$_{\bar \Lambda} =$
 $\{ (p,q): p,q \in {\goth p}_{\bar \Lambda}$ and
 $\[ F \in \tau \wedge$ arity$(F)>1 \Rightarrow p(F) = q(F)\]\}$

 \

 4.6  MAIN  Claim:
   Assume

 (a) {$\bar \Lambda$} is an alphabet sequence for  a vocabulary $\tau =
\tau [{\bar \Lambda} ].$

 (b)  $k_0 \ge  f^6_{{\bar \Lambda}}( \ell+1, c)$  and
$k_0 > \ell.$

 (c)  $K$ is a $\tau-$fim  of dimension $k_0  - 1$
 with $A_2$ the last $\ell$ elements and $A_1$ the first
 $(k_0 -\ell -1) -$elements  (this $K$ surve
just for notation).

 (d)  $\tau^*$ is the vocabulary
  $\tau_{K,A_1,A_2}$, see Definition 1.5(1) and
 proj is the following function
 from $\tau^*$ to $\tau:$ it map $F_{K,\bar a_1, \bar a_2}$
 to  $F$ and {$\bar \Lambda$}$^* =^{df} \langle \Lambda^*_{F} : F \in
\tau^* \rangle$
 where $\Lambda^*_{F} = \Lambda_{{\rm proj}(F)}$,
 so proj $\restriction \tau$ is the identity.

  (e) $c^* =^{df} c^{{\rm card}({\rm Space}_{\bar
\Lambda}(K))}.$

  {\bold Then}

 $f^6_{\bar \Lambda}(\ell,c)  \le k_0 +$
 $f^7_{{\bar \Lambda}^*}
({\Bbb P}_{{\bar \Lambda}^*},1, c^*  )  -1$

 \

REMARK: This is similar to the proof of 2.4,
but for completeness we do it in full.

 Proof:
 Let  $k_1 = f^7_{{\bar \Lambda}^*}(1,
{\Bbb P}_{{\bar \Lambda}^*}, c^*  )$
 and let $k= k_0 + k_1 - 1$, so it suffice to prove
 that $k \ge  f^0_{\bar \Lambda}(n ,\ell , c ).$
 For this it is enough to check (*)$_{k}$
from Definition
4.2, also let
 $M$  be a fim of vocabulary $\tau$
 and dimension  $k$ (that is  $P^{M}$ is
 with  $k$ members), $V =$ Space$_{\bar \Lambda}
(M)$,
 and $d$ an  $(\ell, 1)-$base-invariant
   $C-$colouring   of $V$
 such that $C$ has  $\le  c$ members.
So we just have to prove  that
  the conclusion of Definition 4.2 holds,
that is there is a monochromatic $V-$line.

 Let
 $w_1 =^{df}  \{ a: a\in P^{M}$
 and the number of $b<^{M}$ a is $\ge k_0  - \ell - 1$ but
  is $<k_0 - \ell - 1 +k_1\}$
 hence in $w_1$  there are $k_1$ members,  and let $w_0$ be the set of
 first  $k_0  -\ell -1$  members of $P^{M}$ by  $<^{M}$, lastly
 let $w_2$ be the set of
 the $\ell$ last members of $M$ by $<^{M}.$  So $w_0, w_1, w_2$
 form a convex
 partition of $P^{M}.$

 Now  we let
  $K$ be $M$
 restricted to cl$_{M} (w_0\cup w_2)$,
 (note that this gives no contradiction to the
assumption on
$K$, as concerning $K$ there, only its
 vocabulary and dimension are important and they fit).
  Let   $K^+$ be a fim
 with vocabulary $\tau$  and dimension $k_0$,
 let $g_0 \in$ PHom$(M,K^+)$ be the
 following function from $P^{M}$ onto $P^{K^+}:$
 it maps all the members of $w_1$ to one member of $P^{K^+}$
 which we call $b^*$, it is a one to one order preserving function
 from $w_2$  onto $\{ b \in P^{K^+} :
 b^* <^{K^+} b \}$
 and it is a one to one order preserving function
 from $w_0$  onto
$\{ b \in P^{K^+} : b <^{K^+} b^* \}.$
 Let $g  \in$ Hom$(M,K^+)$ be
the unique extension  of $g_0;$
 without loss of generality $g_0$
is the identity on $w_0$ and on $w_2$
 hence without loss of generality $g$
is the identity on $K$, it exist by 1.2.

 Next recall that the vocabulary $\tau^* = \tau_{K, w_{o},w_2}$
 is a well defined vocabulary (see Definition
1.5(1)).  Next
 we shall define a  $\tau^*-$fim $N.$
 Its universe is  $M \setminus  K ;$
 we let $P^{N}$  be $w_1$  and $<^{N}$  be
 $<^{M} \restriction P^{N}.$
 Now  we have to define the function
$F^{N}_{K,\bar a_1, \bar a_2}$,
 say of arity $r$, where $F \in   \tau, \bar a_1$
 a non decreasing sequence
from $w_0$
 and $\bar a_2$ a non decreasing sequence from $w_2$,
 and lg($\bar a_1)  +$
 lg$(\bar a_2 )  <$ arity$(F)$. So
$r= {\rm arity}(F) - {\rm lg}( \bar a_1 ) -
{\rm lg}( a_2 )$.

 For $b_1 \le^{N} \dots \le^{N} b_{r} \in P^{N}$
 we let  $F^{N}_{\bar a_1, \bar a_2}
( b_1, \dots, b_{r} )$   be
equal to

 $b= F^{M}( \bar a_1, b_1,\dots, b_r, \bar  a_2 )$
 $= F^{M} ( a^1_1, a^1_2,\dots, a^1_{{\rm lg}(\bar a_1 )}$,
 $b_1,\dots, b_{r}, a^2_1, \dots, a^2_{{\rm lg}(\bar a_2) } ).$

 It is easy to check that the
number of arguments is
right
 and also the sequence they form is
$\le^{M}-$increasing, so this is well defined and
belongs to $M$, but still we have to check that it
belongs to $N.$
 But  $N = M \setminus K$ and if $b \in K$ then
  base$_{\lg(\bar a_1 )+1}(b) \in K$
 and it is just $b_1$ which belongs to $w_1$,
contradiction.
 Lastly it is also trivial to note
 that every member of $N$ has this form.
 It is easy to check that $N$ is really a
$\tau^*-$fim.

 We next let  $V^* =$ Space$_{{\bar \Lambda}^*} (N)$ let
$C^* = \{ g: g$ is a function from ${\rm Space}_{\bar \Lambda} (K)$
to  $C\}$  and
 define a $C^*-$ colouring $d^*$ of  $V^*.$
 For $\eta \in V^*$ let $d^*(\eta)$ be the following function
 from Space$_{\bar \Lambda} (K)$ to $C:$
 for $\nu\in K$ we let $\( d^*(\eta) \) (\nu) =
 d(\eta \cup \nu ).$

 Clearly the function $d^*(\eta)$ is a $C^*-$colouring of  $K.$
 How many such functions, that is members of $C^*$ there are?
 The domain  has clearly   card(Space$_{\bar \Lambda}(K))$
 members,
  (we can get slightly less if $\ell >0$,
but with no real influence).
 The range has at most $c$ members, so
the number of such functions
 is at most $c^{{\rm card}({\rm Space}_{\bar \Lambda}(K)) }$,
 a number which we have called $c^*$  in the claim's
statement.

  Hence $d^*$ is a  $c^*-$colouring.

 So as we have chosen $k_1 =$
  $f^7_{{\bar \Lambda}^*}({\Bbb P}_{{\bar \Lambda}^*},1, c^*  )$
 we can apply Definition 4.5
   to $V^* =$ Space$_{{\bar \Lambda}^*}(N)$ and $d^*;$
 so we can find a $d^*-$monochromatic $V^*-$line  $L^*.$
  Let  $h$ be the  function
 from  $U =^{df}$  Space$_{\bar \Lambda} (K^+)$  to $V$ defined as follows:

 (*) $h (\rho)=\nu$ iff :

 (a) $\nu\in V, \rho\in U$,

 (b)  $\nu \restriction    K =\rho \restriction
K$

 (c) if  $b \in  N \setminus$ supp$_{N}(L^* )$
  then $\nu(b) = \eta(b)$ for every
 $\eta \in L^*.$

 (d)  if a $\in$ supp$_{N}(L^*)$, (so a $\in  N$,
 $F_{N,a} = F_{K,\bar a_1,\bar a_2}$,
 base$_{N}(a) \subseteq$ supp$^{P}_{N}(L^*) )$,
  and $b   \in   K^+$,
 $F_{K^+,b} = F, b= F(\bar a_1, b^*,\dots,b^*,\bar a_2)$
 (with  the number of cases of $b^*$ being
 arity$(F_{K, \bar a_1,\bar a_2} ))$
  then
  $\rho(b)=\nu(a).$

 \

 Let the range of $h$ be called $S.$ Now clearly

 $\otimes_1 (\alpha) h$ is a one to one function from  $U$  to $S \subseteq
V.$

 $(\beta)$ $S$ has  $|$Space$_{\bar \Lambda} (K^+ )|$ members

 $(\gamma)$ $S$ is a subspace of $V$ of dimension $k_0,h(\rho) =$
pt$_{S}(\rho)$,
 see  Definition 1.7(7).

 \

 Now clearly

 $\otimes_2$ there is a $C-$colouring $d^{\circ}$ of $U$ such that:

 $d^{\circ} (\nu) = d(h(\nu))$  for $\nu \in U.$

 \

and

 $\otimes_3$  $d^{\circ}$ is $(\ell+1,1)-$base -invariant

 [WHY? Reflect]

 \

 Applying the definition of
$k_0  = f^6_{\bar \Lambda}( \ell +1,c)$ ,
 that is Definition 4.2
 to $\bar \Lambda, U,d^{\circ}$ we can conclude
that there is
 a $d^{\circ} -$monochromatic $U-$line $L^{\circ}.$
 Let
 $L =^{df} \{ h(\rho) : \rho \in L^{\circ} \}.$
 It is easy to check that $L$ is as required.

 $\qed_{4.6}$

 \

 4.7  Claim:
 (1) Assume that $\bar \Lambda$ is a $\tau-$alphabet sequence ,
 and $p^* \in {\goth p}_{\bar \Lambda}$  and {\Bbb P}$^+ =$
 {\Bbb P} $\cup\{(p^*,q):  q \in {\goth p}_{\bar \Lambda}$ and

 $\[ F \in \tau$
 $\wedge$ arity
$(F)>1 \Rightarrow q(F) = p^{*}(F) \} \subseteq$
 {\Bbb P}$_{\bar \Lambda}$ (see Definition 4.5)
 and $n = \Pi_{F \in \tau, {\rm arity}(F) =1}|\Lambda_{F}|.$
 {\bold Then}

 $f^7_{\bar \Lambda}({\Bbb P}^+,m ,c) \le$
 HJ$(n, f^7_{\bar\Lambda} ({\Bbb P},m,c),c)$

 (on HJ see 0.3(2)).

 (2) $f^7(\bar \Lambda,{\Bbb P}_{\bar \Lambda},m,c)$ is in {\Bbb E}$_6.$

 Proof:
 (1) Straight. Let $M$ be a $\tau-$fim of dimension $k =^{df}$
 HJ$(n, f^7_{\bar\Lambda} ({\Bbb P},m,c),c)$  and let $V=$ Space$_{\bar
\Lambda}(M)$
 and let $d$ be a $c-$colouring of $V.$
  Let $\tau^* = \{ F \in \tau :$arity$(F)=1\}$ and let $M^*$ be a $\tau^*-$fim
 of dimension $k$, without loss of generality $M^*$ is $M$ restricted to
$\tau^*$ and the universe of
 $M^*.$ Let  {$\bar \Lambda$}$^* = {\bar \Lambda} \restriction \tau^*$ and let
  $V^* =$ Space$_{{\bar \Lambda}^*}(M^*)$ and let $h$ be the  function from
 $V^*$ to $V$ defined as follows: let $\eta \in V^*$,
 for $b \in M$ let $(h(\eta))(b)$ be $\eta(b)$ is $b \in  M^*$
 and be $p^*(F_{M,b})$ if $b \in  M \setminus M^*.$
So $h$ is a function as
  required and we define  a $c-$colouring $d^*$ of $V^*$ by
 $d^*(\nu)= d(h(\nu))$ for $\nu \in V^*.$

 Now we apply the definition of
$k=$ HJ$(n, f^7_{\bar\Lambda} ({\Bbb P},m,c),c)$
 to the space $V^*$ and the colouring $d^*$
 and we get a subspace $S^*$ of $V^*$
 on which $d^*$ is constant and has dimension
 $f^7_{\bar\Lambda} ({\Bbb P},m,c).$
 There is a unique subspace $S'$ of $V$ of dimension
$f^7_{\bar\Lambda} ({\Bbb P},m,c)$
 such that  $\eta \in S' \Rightarrow \eta
\restriction M^* \in S^*.$
Clearly:

$(*)_1$ if  $L$ is a $V-$line which is $\ps S$
and $(p,q) \in \xP^+ \pb \xP$ then
$pt_L (p) = d ( pt_L (q))$

Now, letting $k' = f^7_{\bar \L} ( \xP , m , c )$
and $ d' =  d \pr S'$, we can apply
the definition of $f^7_{\bar \L} (\xP , m , c )$
and get a subspace $S$ of $S'$ of
dimension $m$ such that

$(*)_2$ if $L$ is a $V-$line which is included in
$S$ and $(p,q) \in \xP$ then
$d' ( pt_L (p)) = d' ( pt_L (q ))$
which means that
$d ( pt_L (p)) = d ( pt_L (q ))$

By $(*)_1 + (*)_2 $, clearly  $S$ is as required.

 (2) Let $\{p^*_{i}: i<i(*)\}$ be maximal
subset of
 {\Bbb P}$_{\bar \Lambda}$ such that
   $ i<j < i(*)
\Rightarrow p^*_{i} \restriction \{F\in \tau:$
arity$(F)>1\} \not=$
 $p^*_{j} \restriction \{F\in \tau:$ arity$(F)>1\}$
  and let {\Bbb P}$_{j} = \{(p^*_{i},q):i<j$  and $q \in {\goth p}_{\bar
\Lambda}$
 and $q\restriction \{F\in \tau:$ arity$(F)>1\}  =$
 $p^*_{i} \restriction \{F\in \tau:$ arity$(F)>1\} \}.$
 By part (1) we have a recursion formula (we use 1.10 freely):

 $f^7_{\bar \Lambda} ({\Bbb P}_{i+1},m,c) \le$
 HJ$(|\Pi_{F\in \tau, {\rm arity}(F)=1}|\Lambda_{F}|,  f^7_{\bar \Lambda}
({\Bbb P}_{i},m,c),c)$

 As HJ belongs to {\Bbb E}$_5$ ( by [Sh 329, 1.8(2),p.691), we are
done.

 $\qed_{4.7}$

 \

 4.8 DEFINITION:  Let $f^{6,*}(\bar \Lambda, \ell,t,c)$
is defined by induction
 on $\ell$ as follows:

 $f^{6,*}(\bar \Lambda,0,t,c) = t$

 $f^{6,*}(\bar \Lambda, \ell+1,t,c)$  is equal to
 $k_0 + f^7_{{\bar \Lambda}[k_0 ]}
({\Bbb P}_{{\bar \Lambda}[k_0 ]} ,
 1,c^{{\rm card}({\rm Space}_{\bar \Lambda}
(M^{\tau}_{k_0})} )  -1$

 where $k_0 = {\rm Max }\{, \ell +1 ,
 f^{6,*}(\bar
\Lambda,\ell,t,c)$
 and {{$\bar \Lambda$}$[k_0 ]$}
 is defined from $\bar \Lambda$ as in the main claim 4.6.

 \

 4.9 Claim: $f^{6,*}$ beongs to {\Bbb E}$_7$

 Proof: Straight.

 \

 4.10  Theorem:
 (1)   The function $f^1({\bar \Lambda},c)$ is well
 defined, i.e. always get value, a natural number and
 is primitive recursive, in fact belongs to {\Bbb E}$_8.$

 (2) Similarly the function $f^6(\bar \Lambda,\ell,c).$

 (3) $f^4$  is primitive recursive, in fact
 belongs to {\Bbb E}$_9.$

 \

 Proof:
 (1),(2) Let $\tau= \tau[{\bar \Lambda}] .$
 The proof follows by induction,
 the main induction is on $t=$  arity$(\tau_{\bar \Lambda} )$
 (or, if you prefer $\Pi_{F \in \tau[{\bar
\Lambda}]}(|\Lambda_{F}|+1) ).$

 CASE 0: arity$(\tau)=1$

 This is Hales-Jewett theorem (on a bound see [Sh:329] or [GRS80])

 CASE 1:  arity$(\tau)>1$

 Let $\tau^*, {\bar \Lambda}^*$ be as in
claim 4.4, so arity$(\tau^*)$
 $\le$ arity$(\tau)/2$  and
$|\tau^*| \le |\tau|\times 2^{{\rm arity}(\tau)}.$

 Let $\ell^* = f^1({\bar \Lambda}^*,c)$  so (by 4.4) clearly
 $f^6_{\bar \Lambda}(\ell^*,c) \le \ell^*$ hence
 (by Definition 4.8) clearly
 $f^{6,*}({\bar \Lambda} , 0,\ell^*,c) = \ell^* =
 f^1({\bar \Lambda}^*,c)$ ;
together we get
$ f^6_{\bar \Lambda}(\ell^*,c) \le
f^{6,*}({\bar \Lambda} , 0,\ell^*,c)$.
 Hence (by 4.6 $+$ Definition 4.8,
we shall prove by induction on
 $\ell \le \ell^*)$ that
 $f^6_{\bar \Lambda}(\ell^* -\ell,c)
\le
f^{6,*}({\bar
\Lambda},\ell,\ell^*,c)$; for $\ell = 0$
this holds by the previous sentnece; for the induction
step, i.e. the proof for $\ell+1$ we apply
Theorem 4.6 with $\ell^* -\ell, \ell^* - (\ell +1)$
here standing for $\ell+1 , \ell$ there
and letting $k_0 = {\rm
Max}\{ \ell^*\ell , f^6_{\bar \L} (\ell^* , c )\}$
and $\t^* , {\bar \L}^* , c^*$ defined as there,
and we get that
$f^6_{\bar \L} ( \ell^* - (\ell +1 ) , c ) \le
k_0 + f^7_{{\bar \L}^*} ( \xP_{{\bar \L}^*} , 1, c^* ) -1
\le  {\rm
Max}\{ \ell^*\ell , f^6_{\bar \L} (\ell^* , c )\} +
f^7_{{\bar \L}^*} ( \xP_{{\bar \L}^*} , 1, c^* ) -1 $

but the last expression is exactly
$f^{6,*} ({\bar \L} , \ell +1 , \ell^* , c )$

 So (using $\ell = \ell^*)$ clearly
 $ f^6_{\bar \Lambda}(0,c)
\le
f^{6,*}({\bar \Lambda},\ell^*,\ell^*,c).$

 Now

 $f^1({\bar \Lambda},c) = f^6_{\bar \Lambda}(0,c) \le$
 $f^{6,*}({\bar \Lambda},\ell^*,\ell^*,c)   \le$
 $f^{6,*}({\bar \Lambda},f^1({\bar \Lambda}^*,c),\varphi^1({\bar
\Lambda}^*,c),c) .$

 As $f^{6,*}$ is from $\xE_7$ by 3.14,
 this clearly give the desired
conclusion.

 (3) Should be clear from the proof of 3.5
and the previous parts.

 $\qed_{4.10}$

\

 REFERENCE

 \

 [BL96]
 V. Bergelson and A. Leibman,
 Polynomial extensions of van der Waerder'$s$
 and Szemeredi theorems,
 JAMS 9(1996)725-753

  [BL 9x] V. Bergelson and A. Leibman,
Set polynomial and polynomial extensions of
the Hales Jewett theorem, to appear

 [Fu81] H. Furstenberg,
 Recurrence in Ergodic Theory and
Combinatorial Number Theory,
 Princeton University Press 1981

 [GR71] R.L. Graham, B.L. Rothschild,
 Ramsey's theorem for $n-$parameter sets,
 TAMS 159(1971)257-292

 [GRS80] R.L. Graham, B.L. Rothschild and
H.J. Spencer,
 Ramsey Theory
 Wiley-Interscience Ser. in Discrete Math.
 New York 1980

 [Ro84] H.E. Rose,
 Subrecursion:  functions and heirarchies,
 Oxford Logic Guide 9,
 Oxford University Press, Oxford 1984

 [Sh:329]
  Shelah, Saharon,
  Primitive recursive bounds for van der Waerden
numbers,
  Journal of the American Mathematical Society,
 1 (1988) 683--697

\end